\documentclass[11pt]{article}
\usepackage{amsthm,amstext,amssymb,amsmath,graphicx,amsfonts,caption,subcaption,mathrsfs,multirow,epigraph,dsfont}
\usepackage[utf8]{inputenc}
\usepackage{booktabs}
\usepackage{tikz}
    \usetikzlibrary{arrows,topaths}
    \usepackage{placeins}
\usepackage[mathscr]{euscript}
\usepackage{multicol}
\usepackage{natbib}
\usepackage[pagebackref=false,colorlinks,linkcolor=blue,citecolor=blue]{hyperref}

\theoremstyle{plain}
\newtheorem{theorem}{Theorem}[section]
\newtheorem{lemma}[theorem]{Lemma}
\newtheorem{proposition}[theorem]{Proposition}

\theoremstyle{definition}
\newtheorem{definition}[theorem]{Definition}
\newtheorem{corollary}[theorem]{Corollary}
\newtheorem{example}[theorem]{Example}

\theoremstyle{remark}
\newtheorem{remark}{\sc Remark}
\makeatother
\makeatletter
\def\namedlabel#1#2{\begingroup
   \def\@currentlabel{#2}%
   \label{#1}\endgroup}
\makeatother
\date{}
\title{\bf  Mp-residuated lattices}\vspace{.25 in}
\author{ \vspace{.25 in} {\bf Saeed Rasouli$^{1*}$} and {\bf Amin Dehghani $^2$}\\
Department of Mathematics, Persian Gulf University, \\Bushehr, Iran \\
{\tt $^1$srasouli@pgu.ac.ir }\\
{\tt $^2$dehghany.amin@hotmail.com}\\}
 
\begin{document}
 \maketitle
%\renewcommand{\epigraphsize}{\large}
% \epigraphnoindent{``\textit{The basic objects of mathematical thought exist only as mental conceptions, though the
%source of these conceptions lies in everyday experience in manifold ways, in the
%processes of counting, ordering, matching, combining, separating, and locating in space
%and time.}"}{\cite{feferman2009conceptions}}
 \begin{abstract}
 This paper is devoted to the study of a fascinating class of residuated lattices, the so-called mp-residuated lattice, in which any prime filter contains a unique minimal prime filter. A combination of algebraic and topological methods is applied to obtain new and structural results on mp-residuated lattices. It is demonstrated that mp-residuated lattices are strongly tied up with the dual hull-kernel topology. Especially, it is shown that a residuated lattice is mp if and only if its minimal prime spectrum, equipped with the dual hull-kernel topology, is Hausdorff if and only if its prime spectrum, equipped with the dual hull-kernel topology, is normal. The class of mp-residuated lattices is characterized by means of pure filters. It is shown that a residuated lattice is mp if and only if its pure filters are precisely its minimal prime filters, if and only if its pure spectrum is homeomorphic to its minimal prime spectrum, equipped with the dual hull-kernel topology.\footnote{2020 Mathematics Subject Classification: 06F99,06D20,06E15 \\
{\it Key words and phrases}: mp-residuated lattice; pure filter; dual hull-kernel topology; pure spectrum.\\
$^*$: Corresponding author}
\end{abstract}
%%%%%%%%%%%%%%%%%%%%%%%%%%%%%%%%%%%%%%%%%%%%%%%%%%%%%%%%%%%%%%%%%%%%%%%%%%%%%%%%%%%%%%%%%
\section{Introduction}

Let $\mathfrak{A}$ be a residuated lattice, $\mathscr{F}(\mathfrak{A})$ the lattice of filters, and $\mathscr{PF}(\mathfrak{A})$ the lattice of principal filters of $\mathfrak{A}$. The lattice of coannihilators of $\mathfrak{A}$, say $\Gamma(\mathfrak{A})$, is the skeleton of $\mathscr{F}(\mathfrak{A})$, and the lattice of coannulets of $\mathfrak{A}$, say $\gamma(\mathfrak{A})$, is the skeleton of $\mathscr{PF}(\mathfrak{A})$.  So $(\Gamma(\mathfrak{A});\vee^{\Gamma},\cap,\{1\},A)$ is a complete Boolean lattice, in which $\vee^{\Gamma}$ is the join in the skeleton, and $\gamma(\mathfrak{A})$ is a sublattice of $\Gamma(\mathfrak{A})$. $\mathfrak{A}$ is said to be \textit{Baer} provided that $\Gamma(\mathfrak{A})$ is a sublattice of $\mathscr{F}(\mathfrak{A})$, and \textit{Rickart} provided that $\gamma(\mathfrak{A})$ is a Boolean sublattice of $\mathscr{F}(\mathfrak{A})$. Obviously, $\mathfrak{A}$ is Rickart if and only if $\gamma(\mathfrak{A})$ is both Boolean and a sublattice of $\mathscr{F}(\mathfrak{A})$. The latter can be characterized by a property that can be formulated in terms of universal algebra, namely that any prime filter contains a unique minimal prime filter.

Historically, this notion is rooted in a query posed by \citet[Problem 70]{birkhoff1940lattice} inspired by M. H. Stone: ``What is the most general pseudocomplemented distributive lattice in which $x^*\vee x^{**}=1$ identically?"  The first solution to this problem belongs to \cite{gratzer1957problem} who gave the name ``\textit{Stone lattices}" to this class of lattices. They characterized stone lattices as distributive pseudocomplemented lattices in which any pair of incomparable minimal prime ideals is comaximal or equivalently each prime ideal contains a unique minimal prime ideal. Motivated by this characterization, \cite{cornish1972normal} studied distributive lattices with zero in which each prime ideal contains a unique minimal prime ideal under the name of ``\textit{normal lattices}". He observed that a distributive lattice with zero, $\mathfrak{A}$, is normal if and only if given $x,y\in A$ such that $x\wedge y=0$, $x^{\perp}$ and $y^{\perp}$ are comaximal. Cornish used this terminology in light of \cite{wallman1938lattices}, who proved that the lattice of closed subsets of a $T_1$ space satisfies the above annihilator condition if and only if the space is normal. \citet[Lemma $\beta$]{artico1976compactness} showed that in a unitary commutative reduced ring any prime ideal contains a unique minimal prime ideal if and only if the set of its annulets is a sublattice of its ideals. \citet[Proposition 2.1]{matlis1983minimal} proved that the class of commutative PF rings, i.e. a unitary ring with the property that every principal ideal is flat, introduced by \citet[p. 151]{hattori1960foundation}, is precisely the class of reduced rings in which any prime ideal contains a unique minimal prime ideal. \citet[Theorem 2.6]{bhattacharjee2013min} tied up the notion of PF rings to the notion of the dual hull-kernel topology. They established that a unitary commutative ring is a PF ring if and only if its minimal prime spectrum, with the dual hull-kernel topology, is Hausdorff. This knot was tightened further by \citet[Theorem 6.2]{aghajani2020characterizations}. They studied the class of unitary commutative rings which fulfill the above universal property, under the name of ``\textit{mp-rings}". They gave a good perspective of mp-rings and asserted that a unitary commutative ring is mp if and only if its prime spectrum, with the dual hull-kernel topology, is normal.

Inspired by the above universal property, many authors have proposed similar notions, under other names, for various structures over the years, see e.g., normal lattices \citep{cignoli1971stone,pawar1994characterizations}, conormal lattices \citep{simmons1980reticulated,belluce1994prime,georgescu2015algebraic}, normal residuated lattices \citep{rasouli2018n}, mp-rings \citep{aghajani2020characterizations}, mp-algebras \citep{georgescu2020}, mp-residuated lattices \citep{rasouli2021mp}, mp-quantales \citep{georgescu2021flat,georgescu2021new}, etc (for a discussion about this terminology, see \citet[p. 185]{simmons1980reticulated} and \citet[p. 78]{johnstone1982stone}).

It is known that residuated lattices play a critical role in the theory of fuzzy logic. Lots of logical algebras such as MTL-algebras, divisible residuated lattices, BL-algebras, MV-algebras, Heyting algebras, and Boolean algebras are subvarieties of residuated lattices. Residuated lattices are not only important from a logical point of view but also interesting from an algebraic point of view and have some interesting algebraic properties.

Given the above discussions, we decided to take a deeper look at mp-residuated lattices. So the notion of  mp-residuated lattices is investigated, and some algebraic and topological characterizations are given. Although, the class of mp-residuated lattices has been investigated by \cite{rasouli2018n}, however, here we give some more characterizations for the class of mp-residuated lattices, which seems to give more light to the topological situation. Our findings show that some results obtained by some above papers can also be reproduced via residuated lattices. Also, outcomes show that mp-residuated lattices can be considered as the dual notion of Gelfand residuated lattices \citep{rasouli2022gelfand}, as asserted in \cite{aghajani2020characterizations} for rings. So mp-residuated lattices can be studied both as one of the two main pillars of Rickart residuated lattices (along with quasicomplemented residuated lattices), and as a dual notion of Gelfand residuated lattices.

This paper is organized into four sections as follows: In Sect. \ref{sec2}, some definitions and facts about residuated lattices are recalled, and some of their propositions extracted. We illustrate this section with some examples of residuated lattices, which will be used in the following sections. Sect. \ref{sec3} deals with mp-residuated lattices. Theorem \ref{nococo} shows that a residuated lattice is mp if and only if the bounded distributive lattice of its filters is conormal. Theorem \ref{noco} (Cornish's characterization) gives an element-wise characterization for mp-residuated lattices. Theorem \ref{osublatfil} shows that a residuated lattice $\mathfrak{A}$ is mp if and only if $\gamma(\mathfrak{A})$ is a sublattice of $\mathscr{F}(\mathfrak{A})$. Theorem \ref{matlchar} (Matlis's characterization) establishes that a residuated lattice $\mathfrak{A}$ is mp if and only if $\mathfrak{A}/D(\mathfrak{p})$ is a domain, for any prime filter $\mathfrak{p}$ of $\mathfrak{A}$. The remaining theorems of this section demonstrate that mp-residuated lattices are strongly tied up with the dual hull-kernel topology. Theorem \ref{mp1culpro} shows that a residuated lattice is mp if and only if its prime spectrum is normal with the dual hull-kernel topology. Sect. \ref{sec4} deals with the pure spectrum of an mp-residuated lattice. The pure filters of an mp-residuated lattice are characterized in Theorem \ref{mppurefcl}. As an important result in this section in Theorem \ref{mp2minspp} is expressed that a residuated lattice is mp if and only if the set of its minimal prime filters is equal to the its purely-prime filters. Theorem \ref{equmpflatmin} verifies that a residuated lattice is mp if and only if the identity map between its pure spectrum and its minimal prime spectrum, equipped with the dual hull-kernel topology, is a homeomorphism. Finally, Corollary \ref{gelspphau} implies that, like Gelfand residuated lattices, the pure spectrum of an mp-residuated lattice is Hausdorff.
%%%%%%%%%%%%%%%%%%%%%%%%%%%%%%%%%%%%%%%%%%%%%%%%%%%%%%%%%%%%%%%%%%%%%%%%%%%%%%%%%%%%%%%%%
\section{Preliminaries}\label{sec2}

In this section, some definitions, properties, and results relative to residuated lattices, which will be used
in the following, recalled.

An algebra $\mathfrak{A}=(A;\vee,\wedge,\odot,\rightarrow,0,1)$ is called a \textit{residuated lattice} provided that $\ell(\mathfrak{A})=(A;\vee,\wedge,0,1)$ is a bounded lattice, $(A;\odot,1)$ is a commutative monoid, and $(\odot,\rightarrow)$ is an adjoint pair. A residuated lattice $\mathfrak{A}$ is called \textit{non-degenerate} if $0\neq 1$. For a residuated lattice $\mathfrak{A}$, and $a\in A$ we put $\neg a:=a\rightarrow 0$ and $a^n:=a\odot\cdots\odot a$ ($n$ times), for any integer $n$. The class of residuated lattices is equational, and so forms a variety. For a survey of residuated lattices, the reader is referred to \cite{galatos2007residuated}.
\begin{remark}\label{resproposition}\citep[Proposition 2.6]{ciungu2006classes}
Let $\mathfrak{A}$ be a residuated lattice. The following conditions are satisfied for any $x,y,z\in A$:
\begin{enumerate}
  \item [$r_{1}$ \namedlabel{res1}{$r_{1}$}] $x\odot (y\vee z)=(x\odot y)\vee (x\odot z)$;
  \item [$r_{2}$ \namedlabel{res2}{$r_{2}$}] $x\vee (y\odot z)\geq (x\vee y)\odot (x\vee z)$.
  \end{enumerate}
\end{remark}
\begin{example}\label{exa6}
Let $A_6=\{0,a,b,c,d,1\}$ be a lattice whose Hasse diagram is given by Figure \ref{figa6}. Routine calculation shows that $\mathfrak{A}_6=(A_6;\vee,\wedge,\odot,\rightarrow,0,1)$ is a residuated lattice in which the commutative operation $``\odot"$ is given by Table \ref{taba6} and the operation $``\rightarrow"$ is given by $x\rightarrow y=\bigvee \{a\in A_6|x\odot a\leq y\}$, for any $x,y\in A_6$.
\FloatBarrier
\begin{table}[ht]
\begin{minipage}[b]{0.56\linewidth}
\centering
\begin{tabular}{ccccccc}
\hline
$\odot$ & 0 & a & b & c & d & 1 \\ \hline
0       & 0 & 0 & 0 & 0 & 0 & 0 \\
        & a & a & a & 0 & a & a \\
        &   & b & a & 0 & a & b \\
        &   &   & c & c & c & c \\
        &   &   &   & d & d & d \\
        &   &   &   &   & 1 & 1 \\ \hline
\end{tabular}
%\begin{tabular}{@{}lllllll@{}}
%\multicolumn{1}{l|}{$\odot$} & 0                      & a                      & b                      & c                      & d                      & 1                      \\ \midrule
%\multicolumn{1}{l|}{0}       & 0                      & 0                      & 0                      & 0                      & 0                      & \multicolumn{1}{l|}{0} \\ \cmidrule(lr){2-2}
%                             & \multicolumn{1}{l|}{a} & a                      & a                      & 0                      & a                      & \multicolumn{1}{l|}{a} \\ \cmidrule(lr){3-3}
%                             &                        & \multicolumn{1}{l|}{b} & a                      & 0                      & a                      & \multicolumn{1}{l|}{b} \\ \cmidrule(lr){4-4}
%                             &                        &                        & \multicolumn{1}{l|}{c} & c                      & c                      & \multicolumn{1}{l|}{c} \\ \cmidrule(lr){5-5}
%                             &                        &                        &                        & \multicolumn{1}{l|}{d} & d                      & \multicolumn{1}{l|}{d} \\ \cmidrule(lr){6-6}
%                             &                        &                        &                        &                        & \multicolumn{1}{l|}{1} & \multicolumn{1}{l|}{1} \\ \cmidrule(l){7-7}
%\end{tabular}
\caption{Cayley table for ``$\odot$" of $\mathfrak{A}_6$}
\label{taba6}
\end{minipage}\hfill
\begin{minipage}[b]{0.6\linewidth}
\centering
  \begin{tikzpicture}[>=stealth',semithick,auto]
    \tikzstyle{subj} = [circle, minimum width=6pt, fill, inner sep=0pt]
    \tikzstyle{obj}  = [circle, minimum width=6pt, draw, inner sep=0pt]

    \tikzstyle{every label}=[font=\bfseries]

    % Before diagram .........................
    \node[subj,  label=below:0] (0) at (0,0) {};
    \node[subj,  label=below:c] (c) at (-1,1) {};
    \node[subj,  label=below:a] (a) at (1,.5) {};
    \node[subj,  label=below right:b] (b) at (1,1.5) {};
    \node[subj,  label=below:d] (d) at (0,2) {};
    \node[subj,  label=below right:1] (1) at (0,3) {};

    \path[-]   (0)    edge                node{}      (a);
    \path[-]   (a)    edge                node{}      (b);
    \path[-]   (0)    edge                node{}      (c);
    \path[-]   (c)    edge                node{}      (d);
    \path[-]   (b)    edge                node{}      (d);
    \path[-]   (d)    edge                node{}      (1);
\end{tikzpicture}

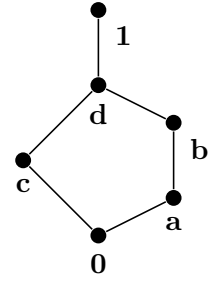
\captionof{figure}{Hasse diagram of $\mathfrak{A}_{6}$}
\label{figa6}
\end{minipage}
\end{table}
\FloatBarrier
\end{example}
\begin{example}\label{exa8}
Let $A_8=\{0,a,b,c,d,e,f,1\}$ be a lattice whose Hasse diagram is given by Figure \ref{figa8}. Routine calculation shows that $\mathfrak{A}_8=(A_8;\vee,\wedge,\odot,\rightarrow,0,1)$ is a residuated lattice in which the commutative operation $``\odot"$ is given by Table \ref{taba8} and the operation $``\rightarrow"$ is given by $x\rightarrow y=\bigvee \{a\in A_8|x\odot a\leq y\}$, for any $x,y\in A_8$.
\FloatBarrier
\begin{table}[ht]
\begin{minipage}[b]{0.56\linewidth}
\centering
\begin{tabular}{ccccccccc}
\hline
$\odot$ & 0 & a & b & c & d & e & f & 1 \\ \hline
0       & 0 & 0 & 0 & 0 & 0 & 0 & 0 & 0 \\
        & a & a & 0 & a & a & a & a & a \\
        &   & b & 0 & 0 & 0 & 0 & b & b \\
        &   &   & c & c & a & c & a & c \\
        &   &   &   & d & a & a & d & d \\
        &   &   &   &   & e & c & d & e \\
        &   &   &   &   &   & f & f & f \\
        &   &   &   &   &   &   & 1 & 1 \\ \hline
\end{tabular}
\caption{Cayley table for ``$\odot$" of $\mathfrak{A}_8$}
\label{taba8}
\end{minipage}\hfill
\begin{minipage}[b]{0.6\linewidth}
\centering
  \begin{tikzpicture}[>=stealth',semithick,auto]
    \tikzstyle{subj} = [circle, minimum width=6pt, fill, inner sep=0pt]
    \tikzstyle{obj}  = [circle, minimum width=6pt, draw, inner sep=0pt]

    \tikzstyle{every label}=[font=\bfseries]

    % Before diagram .........................
    \node[subj,  label=below:0] (0) at (0,0) {};
    \node[subj,  label=below:a] (a) at (-1,1) {};
    \node[subj,  label=below:b] (b) at (1,1) {};
    \node[subj,  label=below:c] (c) at (-2,2) {};
    \node[subj,  label=below:d] (d) at (0,2) {};
    \node[subj,  label=below:e] (e) at (-1,3) {};
    \node[subj,  label=below:f] (f) at (1,3) {};
    \node[subj,  label=below:1] (1) at (0,4) {};

    \path[-]   (0)    edge                node{}      (a);
    \path[-]   (0)    edge                node{}      (b);
    \path[-]   (b)    edge                node{}      (d);
    \path[-]   (d)    edge                node{}      (f);
    \path[-]   (f)    edge                node{}      (1);
    \path[-]   (a)    edge                node{}      (d);
    \path[-]   (a)    edge                node{}      (c);
    \path[-]   (c)    edge                node{}      (e);
    \path[-]   (d)    edge                node{}      (e);
    \path[-]   (e)    edge                node{}      (1);
\end{tikzpicture}

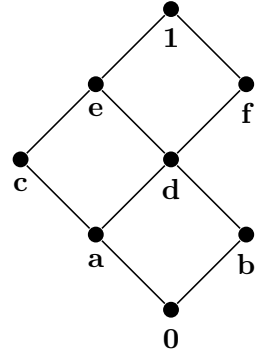
\captionof{figure}{Hasse diagram of $\mathfrak{A}_{8}$}
\label{figa8}
\end{minipage}
\end{table}
\FloatBarrier
\end{example}

Let $\mathfrak{A}$ be a residuated lattice. A non-void subset $F$ of $A$ is called a \textit{filter} of $\mathfrak{A}$ provided that $x,y\in F$ implies $x\odot y\in F$, and $x\vee y\in F$, for any $x\in F$ and $y\in A$. The set of filters of $\mathfrak{A}$ is denoted by $\mathscr{F}(\mathfrak{A})$. A filter $F$ of $\mathfrak{A}$ is called \textit{proper} if $F\neq A$. For any subset $X$ of $A$, the \textit{filter of $\mathfrak{A}$ generated by $X$} is denoted by $\mathscr{F}(X)$. For each $x\in A$, the filter generated by $\{x\}$ is denoted by $\mathscr{F}(x)$ and said to be \textit{principal}. The set of principal filters is denoted by $\mathscr{PF}(\mathfrak{A})$. Following \citet[\S 5.7]{gratzer2011lattice}, a join-complete lattice $\mathfrak{A}$, is called a \textit{frame} if it satisfies the join infinite distributive law (JID), i.e., for any $a\in A$ and $S\subseteq A$, $a\wedge \bigvee S=\bigvee \{a\wedge s\mid s\in S\}$. A frame $\mathfrak{A}$ is called complete provided that $\mathfrak{A}$ is a complete lattice. According to \cite{galatos2007residuated}, $(\mathscr{F}(\mathfrak{A});\cap,\veebar,\textbf{1},A)$ is a complete frame, in which $\veebar \mathcal{F}=\mathscr{F}(\cup \mathcal{F})$, for any $\mathcal{F}\subseteq \mathscr{F}(\mathfrak{A})$.
\begin{example}\label{filterexa}
Consider the residuated lattice $\mathfrak{A}_6$ from Example \ref{exa6} and the residuated lattice $\mathfrak{A}_8$ from Example \ref{exa8}. The sets of their filters are presented in Table \ref{tafiex}.
\begin{table}[h]
\centering
\begin{tabular}{ccl}
\hline
                 & \multicolumn{2}{c}{Filters}                                       \\ \hline
$\mathfrak{A}_6$ & \multicolumn{2}{c}{$\{1\},\{a,b,d,1\},\{c,d,1\},\{d,1\},A_6$} \\
$\mathfrak{A}_8$ & \multicolumn{2}{c}{$\{1\},\{a,c,d,e,f,1\},\{c,e,1\},\{f,1\},A_8$} \\ \hline
\end{tabular}
\caption{The sets of filters of $\mathfrak{A}_6$ and $\mathfrak{A}_8$}
\label{tafiex}
\end{table}
\end{example}

The proof of the following proposition has a routine verification, and so it is left to the reader.
\begin{proposition}\label{genfilprop}
Let $\mathfrak{A}$ be a residuated lattice and $F$ be a filter of $\mathfrak{A}$. The following assertions hold, for any $x,y\in A$:
\begin{enumerate}
  \item  [(1) \namedlabel{genfilprop1}{(1)}] $\mathscr{F}(x)=\{a\in A|x^n\leq a,\textrm{~for~some~integer~}n\}$;
  \item  [(2) \namedlabel{genfilprop2}{(2)}] $x\leq y$ implies $\mathscr{F}(y)\subseteq \mathscr{F}(x)$.
  \item  [(3) \namedlabel{genfilprop3}{(3)}] $\mathscr{F}(x)\cap \mathscr{F}(y)=\mathscr{F}(x\vee y)$;
  \item  [(4) \namedlabel{genfilprop4}{(4)}] $\mathscr{F}(x)\veebar \mathscr{F}(y)=\mathscr{F}(x\odot y)$;
  \item  [(5) \namedlabel{genfilprop5}{(5)}] $\mathscr{PF}(\mathfrak{A})$ is a sublattice of $\mathscr{F}(\mathfrak{A})$.
  %\item  [(6) \namedlabel{genfilprop6}{(6)}] if $\mathfrak{A}$ is finite, then $\mathscr{F}(\mathfrak{A})=\mathscr{PF}(\mathfrak{A})$.
\end{enumerate}
\end{proposition}

The following proposition gives a characterization for the comaximal filters of a residuated lattice.
\begin{proposition}\label{compropo}
  Let $\mathfrak{A}$ be a residuated lattice and $F,G$ two proper filters of $\mathfrak{A}$. The following assertions are equivalent:
  \begin{enumerate}
  \item  [$(1)$ \namedlabel{compropo1}{$(1)$}] $F$ and $G$ are comaximal, i.e., $F\veebar G=A$;
  \item  [$(2)$ \namedlabel{compropo2}{$(2)$}] there exist $f\in F$ and $g\in G$ such that $f\odot g=0$;
  \item  [$(3)$ \namedlabel{compropo3}{$(3)$}] there exists $a\in A$ such that $a\in F$ and $\neg a\in G$.
\end{enumerate}
\end{proposition}
\begin{proof}
\item [\ref{compropo1}$\Rightarrow$\ref{compropo2}:] It is evident by Proposition \ref{genfilprop}.
\item [\ref{compropo2}$\Rightarrow$\ref{compropo3}:] Let $f\odot g=0$, for some $f\in F$ and $g\in G$. This implies that $g\leq \neg f$, and this holds the result.
\item [\ref{compropo3}$\Rightarrow$\ref{compropo1}:] It is evident.
\end{proof}

Let $\mathfrak{A}$ be a residuated lattice. A maximal element in the set of proper filters of $\mathfrak{A}$ is called \textit{maximal}, and the set of maximal filters of $\mathfrak{A}$ denoted by $Max(\mathfrak{A})$. A  meet-irreducible element in the set of proper filters of $\mathfrak{A}$ is called \textit{prime}, and the set of prime filters of $\mathfrak{A}$ denoted by $Spec(\mathfrak{A})$. Since $\mathscr{F}(\mathfrak{A})$ is a distributive lattice, so $Max(\mathfrak{A})\subseteq Spec(\mathfrak{A})$. Zorn's lemma verifies that any proper filter is contained in a maximal filter, and so in a prime filter.

A non-empty subset $\mathscr{C}$ of $\mathfrak{A}$ is called \textit{$\vee$-closed} if it is closed under the join operation, i.e $x,y\in \mathscr{C}$ implies $x\vee y\in \mathscr{C}$.
\begin{theorem}\cite[Theorem 3.18]{rasouli2019going}\label{prfilth}
If $\mathscr{C}$ is a $\vee$-closed subset of $\mathfrak{A}$ which does not meet the filter $F$, then $F$ is contained in a filter $P$ which is maximal with respect to the property of not meeting $\mathscr{C}$; furthermore $P$ is prime.
\end{theorem}

 A minimal element in the set of prime filters of a residuated lattice $\mathfrak{A}$ is called \textit{minimal prime}, and the set of minimal prime filters of $\mathfrak{A}$ denoted by $Min(\mathfrak{A})$. For the basic facts concerning prime filters of a residuated lattice, the reader is referred to \cite{rasouli2019going}.
\begin{example}\label{maxminex}
Consider the residuated lattice $\mathfrak{A}_6$ from Example \ref{exa6} and the residuated lattice $\mathfrak{A}_8$ from Example \ref{exa8}. The sets of their maximal, prime, and minimal prime filters are presented in Table \ref{prfiltab}.
\begin{table}[h]
\centering
\begin{tabular}{cccc}
\hline
                 & \multicolumn{3}{c}{Prime filters}      \\ \hline
                 & Maximal filters &          & Minimal prime filters      \\
$\mathfrak{A}_6$ &  $\{a,b,d,1\},\{c,d,1\}$                 & $\{d,1\}$ &$ \{1\}$\\
$\mathfrak{A}_8$ &  $\{a,c,d,e,f,1\}$                          &         & $\{c,e,1\},\{f,1\}$ \\ \hline
\end{tabular}
\caption{The sets of maximal, prime, and minimal prime filters of $\mathfrak{A}_6$ and $\mathfrak{A}_8$}
\label{prfiltab}
\end{table}
\end{example}
\begin{proposition}\citep{rasouli2019going}\label{1mineq}
Let $\mathfrak{A}$ be a residuated lattice. The following assertions hold:
\begin{enumerate}
  \item  [$(1)$ \namedlabel{1mineq1}{$(1)$}] A subset $P$ of $A$ is a minimal prime filter if and only if $\dot{P}\stackrel{def.}{=}A\setminus P$ is a $\vee$-closed subset of $\mathfrak{A}$ which it is maximal with respect to the property of not containing $1$;
  \item  [$(2)$ \namedlabel{1mineq2}{$(2)$}] any prime filter of a residuated lattice contains a minimal prime filter;
  \item  [$(3)$ \namedlabel{1mineq3}{$(3)$}] a prime filter $P$ of $\mathfrak{A}$ is minimal prime if and only if  for any $x\in A$, $P$ contains precisely one of $x$ or $x^{\perp}$.
\end{enumerate}
\end{proposition}

Let $\mathfrak{A}$ be a residuated lattice and $\Pi$ a collection of prime filters of $\mathfrak{A}$.
For a subset $\pi$ of $\Pi$ we set $k(\pi)=\bigcap \pi$, and for a subset $X$ of $A$ we set $h_{\Pi}(X)=\{P\in \Pi\mid X\subseteq  P\}$ and $d_{\Pi}(X)=\Pi\setminus h_{\Pi}(X)$. The collection $\Pi$ can be topologized by taking the collection $\{h_{\Pi}(x)\mid x\in A\}$  as a closed (an open) basis, which is called \textit{the (dual) hull-kernel topology} on $\Pi$ and denoted by $\Pi_{h(d)}$. The generated
topology by $\tau_{h}\cup \tau_{d}$ on $Spec(\mathfrak{A})$ is called \textit{the patch topology} and denoted by $\tau_{p}$. As usual, the Boolean lattice of all clopen subsets of a topological space $A_{\tau}$ shall be denoted by $Clop(A_{\tau})$. For a detailed discussion on the (dual) hull-kernel and patch topologies on a residuated lattice, we refer to \cite{rasouli2018hull}.
\begin{proposition}\citep{rasouli2018hull}\label{minpro}
Let $\mathfrak{A}$ be a residuated lattice. We have:
\[Clop(Spec_{d}(\mathfrak{A}))=\{h(e)\mid e\in \beta(\mathfrak{A})\}.\]
\end{proposition}

Let $\Pi$ be a collection of prime filters in a residuated lattice $\mathfrak{A}$. In the following, for a given subset $\pi$ of $\Pi$, $cl^{\Pi}_{h(d)}(\pi)$ stands for the closure of $\pi$ in the topological space $(\Pi,\tau_{h(d)})$. If $\pi=\{P\}$ for some prime filter $P$ of $\mathfrak{A}$, then $cl^{\Pi}_{h(d)}(\{P\})$ is simply denoted by $cl^{\Pi}_{h(d)}(P)$. If $\Pi$ is understood, it will be dropped.
\begin{lemma}\citep[Theorem 3.14]{rasouli2018hull}\label{retractlemma}
Let $\mathfrak{A}$ be a residuated lattice, $\Pi$ a collection of prime filters of $\mathfrak{A}$ and $\mathfrak{p},\mathfrak{q}\in \Pi$. The following assertions are equivalent:
\begin{enumerate}
   \item [(1) \namedlabel{retractlemma1}{(1)}] $\mathfrak{p}\subseteq \mathfrak{q}$;
   \item [(2) \namedlabel{retractlemma2}{(2)}] $\mathfrak{q}\in cl_{h}(\mathfrak{p})$;
   \item [(3) \namedlabel{retractlemma3}{(3)}] $\mathfrak{p}\in cl_{d}(\mathfrak{q})$.
 \end{enumerate}
\end{lemma}

The following proposition characterizes the open sets of the spectrum of a residuated lattice w.r.t the dual hull-kernel topology.
\begin{proposition}\label{opensd}
  Let $\mathfrak{A}$ be a residuated lattice. The open sets of $Spec_{d}(\mathfrak{A})$ are precisely of the form $\{\mathfrak{p}\in Spec(\mathfrak{A})\mid \mathfrak{p}\cap X\neq \emptyset\}$, where $X$ is a subset of $A$.
\end{proposition}
\begin{proof}
  Let $U$ be an open set in $Spec_{d}(\mathfrak{A})$. So $U=\bigcup_{x\in X} h(x)$, for some $X\subseteq A$. It is clear that $\bigcup_{x\in X} h(x)=\{\mathfrak{p}\in Spec(\mathfrak{A})\mid \mathfrak{p}\cap X\neq \emptyset\}$. This holds the result.
\end{proof}
\begin{remark}\label{closd}
  Let $\mathfrak{A}$ be a residuated lattice. By Proposition \ref{opensd} follows that the closed sets of $Spec_{d}(\mathfrak{A})$ are precisely of the form $\{\mathfrak{p}\in Spec(\mathfrak{A})\mid \mathfrak{p}\cap X=\emptyset\}$, where $X$ is a subset of $A$.
\end{remark}

Let $\Pi$ be a collection of prime filters in a residuated lattice $\mathfrak{A}$. Following \citet[p. 290]{de1983projectivity}, if $\pi$ is a subset of $\Pi$, its \textit{specialization} (\textit{generalization}) in $\Pi$, $\mathscr{S}_{\Pi}(\pi)$ ($\mathscr{G}_{\Pi}(\pi)$), is the set of all primes in $\Pi$, which contain (are contained in) some prime belonging to $\pi$. One can see that $\mathscr{S}$ and $\mathscr{G}$ are closure operators on the power set of $Spec(\mathfrak{A})$. A fixed point of $\mathscr{S(G)}$ is called $\mathscr{S}_{\Pi}$-stable ($\mathscr{G}_{\Pi}$-stable). If $\Pi$ is understood, it will be dropped. Notice that for any subset $B$ of $A$, $\bigcup_{b\in B}h(b)(\bigcup_{b\in B}d(b))$ is $\mathscr{S(G)}$-stable. The following theorem characterizes the closed sets of the (dual) hull-kernel topology.
\begin{theorem}\cite[Theorem 4.30]{rasouli2018hull}\label{closefalzai}
  Let $\mathfrak{A}$ be a residuated lattice and $\pi$ a subset of $Spec(\mathfrak{A})$. $\pi$ is closed under the dual hull-kernel topology if and only if it is closed under the patch topology and $\mathscr{G}$-stable.
\end{theorem}

For a residuated lattice $\mathfrak{A}$ the hull-kernel topology on $Min(\mathfrak{A})$ is a well-studied structure. For example, it is known that the hull-kernel topology on $Min(\mathfrak{A})$ is totally disconnected \cite[Corollary 5.5]{rasouli2018hull}, and classifications of when $Min(\mathfrak{A})$ is compact \cite[Theorem 5.10]{rasouli2018hull}. In the sequel, we fucose on the dual hull-kernel topology on $Min(\mathfrak{A})$. In particular, we characterize when $Min_{d}(\mathfrak{A})$ is Hausdorff.
\begin{proposition}\citep[Theorem 4.6(2)]{rasouli2018hull}\label{dhulkerinstr}
Let $\mathfrak{A}$ be a residuated lattice. $Spec_{d}(\mathfrak{A})$ and $Min_{d}(\mathfrak{A})$ are compact.
\end{proposition}

Let $\mathfrak{A}$ be a residuated lattice. For any subset $X$ of $A$, we set $X^{\perp}=kd(X)$, $\Gamma(\mathfrak{A})=\{X^{\perp}\mid X\subseteq A\}$, $\gamma(\mathfrak{A})=\{x^{\perp}\mid x\in A\}$, and $\lambda(\mathfrak{A})=\{x^{\perp\perp}\mid x\in A\}$. Elements of $\Gamma(\mathfrak{A})$, $\gamma(\mathfrak{A})$, and $\lambda(\mathfrak{A})$ are called \textit{coannihilators}, \textit{coannulets}, and \textit{dual coannulets} of $\mathfrak{A}$, respectively.

Let $\mathfrak{A}$ be a $\wedge$-semilattice with zero. Recall \citep[\S I.6.2]{gratzer2011lattice} that an element $a^{*}\in A$ is a \textit{pseudocomplement} of $a\in A$ if $a\wedge a^{*}=0$ and $a\wedge x=0$ implies that $x\leq a^{*}$. An element can have at most one pseudocomplement. $\mathfrak{A}$ is called \textit{pseudocomplemented} if every element of $A$ has a pseudocomplement. The set $S(\mathfrak{A})=\{a^{*}\mid a\in A\}$ is called \textit{the skeleton of} $\mathfrak{A}$ and we have $S(\mathfrak{A})=\{a\in A\mid a=a^{**}\}$. By \citet[Theorem 100]{gratzer2011lattice} follows that if $\mathfrak{A}$ is a pseudocomplemented complete $\wedge$-semilattice, then $S(\mathfrak{A})$ is a complete Boolean lattice, where the meet in $S(\mathfrak{A})$ is calculated in $\mathfrak{A}$, the join in $S(\mathfrak{A})$ is given by $\vee X=(\wedge \{x^{*}\mid x\in X\})^{*}$, for any $X\subseteq S(\mathfrak{A})$, and $1\stackrel{def.}{=}0^*$.

Applying Proposition 2.11 of \cite{rasouli2020quasicomplemented}, it follows that $\Gamma(\mathfrak{A})$ is the skeleton of $\mathscr{F}(\mathfrak{A})$, and $\gamma(\mathfrak{A})$ is the skeleton of $\mathscr{PF}(\mathfrak{A})$. So $(\Gamma(\mathfrak{A});\vee^{\Gamma},\cap,\{1\},A)$ is a complete Boolean lattice, in which $\vee^{\Gamma}$ is the join in the skeleton, and $\gamma(\mathfrak{A})$ is a sublattice of $\Gamma(\mathfrak{A})$. $\mathfrak{A}$ is said to be \textit{Baer} provided that $\Gamma(\mathfrak{A})$ is a sublattice of $\mathscr{F}(\mathfrak{A})$, and \textit{Rickart} provided that $\gamma(\mathfrak{A})$ is a Boolean sublattice of $\mathscr{F}(\mathfrak{A})$. For the basic facts concerning coannihilators and coannulets of residuated lattices we refer to \cite{rasouli2018generalized}.

   Let $\mathfrak{A}$ be a residuated lattice. For a $\vee$-closed subset $I$ of $\ell(\mathfrak{A})$, set $\omega(I)=\{a\in A|a\vee x=1,\textrm{~for~some~} x\in I\}$, and $\Omega(\mathfrak{A})=\{\omega(I)|I\in id(\ell(\mathfrak{A}))\}$. Using Proposition 3.4 of \cite{rasouli2018n}, it follows that $\Omega(\mathfrak{A})\subseteq \mathscr{F}(\mathfrak{A})$, and so elements of $\Omega(\mathfrak{A})$ are called \textit{$\omega$-filters} of $\mathfrak{A}$. For an $\omega$-filter $F$ of $\mathfrak{A}$, $I_{F}$ denoted an ideal of $\ell(\mathfrak{A})$, which satisfies $F=\omega(I_{F})$. \citet[Proposition 3.7]{rasouli2018n} show that $(\Omega(\mathfrak{A});\cap,\vee^{\omega},\{1\},A)$ is a bounded distributive lattice, in which $F\vee^{\omega} G=\omega(I_{F}\curlyvee I_{G})$, for any $F,G\in \Omega(\mathfrak{A})$ (by $\curlyvee$, we mean the join operation in the lattice of ideals of $\ell(\mathfrak{A})$). For any proper filter $H$ of $\mathfrak{A}$ we set $D(H)=\omega(\dot{H})$. Elements of $D(\{1\})$ shall be called the \textit{unit divisors of $\mathfrak{A}$}. For the basic facts concerning $\omega$-filters of a residuated lattice, interested readers are referred to \cite{rasouli2018n}.
\begin{proposition}\citep{rasouli2018n}\label{omegprop}
    Let $\mathfrak{A}$ be residuated lattice. The following assertions hold:
    \begin{enumerate}
\item  [$(1)$ \namedlabel{omegprop1}{$(1)$}] $\gamma(\mathfrak{A})$ is a sublattice of $\Omega(\mathfrak{A})$;%  \cite[Proposition 3.9]{rasouli2018n};
\item  [$(2)$ \namedlabel{omegprop2}{$(2)$}]  $D(\mathfrak{p})=k\mathscr{G}(\mathfrak{p})=k(\mathscr{G}(\mathfrak{p})\cap Min(\mathfrak{A}))$, for any prime filter $\mathfrak{p}$ of $\mathfrak{A}$;% \cite[Corollary 3.23]{rasouli2018n};
\item  [$(3)$ \namedlabel{omegprop3}{$(3)$}] a prime filter $\mathfrak{p}$ of $\mathfrak{A}$ is minimal prime if and only if $\mathfrak{p}=D(\mathfrak{p})$.
%\item  [$(4)$ \namedlabel{omegprop4}{$(4)$}] $\Omega(\mathfrak{A})\subseteq \alpha(\mathfrak{A})$; %\cite[Proposition 5.11]{rasouli2020quasicomplemented};
%\item  [$(5)$ \namedlabel{omegprop5}{$(5)$}] $\omega(X)=A$ if and only if $1\in X$.% \cite[Proposition 3.2(6)]{rasouli2018n};
%\item  [$(6)$ \namedlabel{omegprop6}{$(6)$}] $D(P)=
%\item  [$(7)$ \namedlabel{omegprop7}{$(7)$}]
\end{enumerate}
  \end{proposition}

\begin{definition}
  A residuated lattice $\mathfrak{A}$ is said to be \textit{a domain} provided that it has no unit divisors.
\end{definition}

The following proposition has a routine verification, and so its proof is left to the reader.
\begin{proposition}\label{domchar}
Let $\mathfrak{A}$ be a residuated lattice and $F$ a filter of $\mathfrak{A}$. The quotient residuated lattice $\mathfrak{A}/F$ is a domain if and only if $F$ is prime.
\end{proposition}
%%%%%%%%%%%%%%%%%%%%%%%%%%%%%%%%%%%%%%%%%%%%%%%%%%%%%%%%%%%%%%%%%%%%%%%%%%%%%%%%%%%%%%%%%%%%%%%%%%%%%%%%%%%%%%
\section{Mp-residuated lattices}\label{sec3}

In this section, the notion of an mp-residuated lattice is investigated, and some topological characterizations of them are extracted.
\begin{definition}\label{nordef}
A residuated lattice $\mathfrak{A}$ is called \textit{mp} provided that any prime filter of $\mathfrak{A}$ contains a unique minimal prime filter of $\mathfrak{A}$.
\end{definition}
\begin{example}\label{quanorexas}
One can see that the residuated lattice $\mathfrak{A}_6$ from Example \ref{exa6} is mp and the residuated lattice $\mathfrak{A}_8$ from Example \ref{exa8} is not mp.
\end{example}
\begin{example}
   The class of MTL-algebras, and so,  MV-algebras, BL-algebras, and Boolean algebras are some subclasses of mp-residuated lattices.
\end{example}

Let $\mathfrak{A}$ be a bounded distributive lattice. $\mathfrak{A}$ is said to be:
\begin{itemize}
  \item \textit{normal} provided that for all $x,y\in A$, $x\vee y=1$ implies there exist $u,v\in A$ such that $u\vee x=v\vee y=1$ and $u\wedge v=0$;
  \item \textit{conormal} provided that for all $x,y\in A$, $x\wedge y=0$ implies there exist $u,v\in L$ such that $u\wedge x=v\wedge y=0$ and $u\vee v=1$.
\end{itemize}
\begin{remark}
In \cite{cornish1972normal} and \cite{pawar1994characterizations}, the above nomenclatures are reversed. We have picked the version of these definitions from \citet[Definition 4.3]{simmons1980reticulated} and \citet[p. 67]{johnstone1982stone} because of the author's discussion in \cite[p. 78]{johnstone1982stone}.
\end{remark}

The following result shows that a residuated lattice is mp if and only if the bounded distributive lattice of its filters is conormal.
\begin{proposition}\label{nococo}
Let $\mathfrak{A}$ be a residuated lattice. The following assertions are equivalent:
\begin{enumerate}
\item  [$(1)$ \namedlabel{nococo1}{$(1)$}] The bounded distributive lattice $\mathscr{F}(\mathfrak{A})$ is conormal;
\item  [$(2)$ \namedlabel{nococo2}{$(2)$}] the bounded distributive lattice $\mathscr{PF}(\mathfrak{A})$ is conormal;
\item  [$(3)$ \namedlabel{nococo3}{$(3)$}] $\mathfrak{A}$ is mp.
\end{enumerate}
\end{proposition}
\begin{proof}
\item [\ref{nococo1}$\Rightarrow$\ref{nococo2}:] Let $x,y\in A$, such that $\mathscr{F}(x)\cap \mathscr{F}(y)=\{1\}$. Then there exist $F,G\in \mathscr{F}(\mathfrak{A})$ such that $F\veebar G=A$ and $F\cap \mathscr{F}(x)=G\cap \mathscr{F}(y)=\{1\}$. Thus there exist $f\in F$ and $g\in G$ such that $f\odot g=0$. This implies that $\mathscr{F}(f)\veebar \mathscr{F}(g)=A$ and $\mathscr{F}(f)\cap \mathscr{F}(x)=\mathscr{F}(g)\cap \mathscr{F}(y)=\{1\}$.
\item [\ref{nococo2}$\Rightarrow$\ref{nococo3}:] Using Proposition \ref{genfilprop}, it is straightforward.
%Let $x,y\in A$, such that $x\vee y=1$. So $\mathscr{F}(x)\cap \mathscr{F}(y)=\{1\}$, and this implies that there exist $f,g\in A$ such that $\mathscr{F}(f)\veebar \mathscr{F}(g)=A$ and $\mathscr{F}(f)\cap \mathscr{F}(x)=\mathscr{F}(g)\cap \mathscr{F}(y)=\{1\}$.
\item [\ref{nococo3}$\Rightarrow$\ref{nococo1}:] Let $F$ and $G$ be two filters of $\mathfrak{A}$ such that $F\cap G=\{1\}$. By distributivity of $\mathfrak{A}$, with a little bit of effort, we can show that $F^{\perp}\veebar G^{\perp}=A$.
\end{proof}

The following theorem gives some algebraic criteria for mp-residuated lattices, inspired by the
one obtained for normal lattices \cite[Theorem 2.4]{cornish1972normal}.
\begin{theorem}\label{noco}(Cornish's characterization)
Let $\mathfrak{A}$ be a residuated lattice. The following assertions are equivalent:
\begin{enumerate}
\item  [$(1)$ \namedlabel{noco1}{$(1)$}] Any two distinct minimal prime filters are comaximal;
\item  [$(2)$ \namedlabel{noco2}{$(2)$}] $\mathfrak{A}$ is mp;
\item  [$(3)$ \namedlabel{noco3}{$(3)$}] for any prime filter $\mathfrak{p}$ of $\mathfrak{A}$, $D(\mathfrak{p})$ is a prime filter of $\mathfrak{A}$;
\item  [$(4)$ \namedlabel{noco4}{$(4)$}] for any maximal filter $\mathfrak{m}$ of $\mathfrak{A}$, $D(\mathfrak{m})$ is a prime filter of $\mathfrak{A}$;
\item  [$(5)$ \namedlabel{noco5}{$(5)$}] for any pairwise elements $x$ and $y$ in $A$, i.e, $x\vee y=1$, $x^{\perp}\veebar y^{\perp}=A$;
\item  [$(6)$ \namedlabel{noco6}{$(6)$}] for any pairwise elements $x$ and $y$ in $A$, there exists $a\in A$ such that $a\in x^{\perp}$ and $\neg a\in y^{\perp}$;
\item  [$(7)$ \namedlabel{noco7}{$(7)$}] for any $x,y\in A$, $(x\vee y)^{\perp}=x^{\perp}\veebar y^{\perp}$;
\item  [$(8)$ \namedlabel{noco8}{$(8)$}] for any $x,y\in A$, $(x\vee y)^{\perp}=A$ implies $x^{\perp}\veebar y^{\perp}=A$.
\end{enumerate}
\end{theorem}
\begin{proof}
\item [\ref{noco1}$\Rightarrow$\ref{noco2}:] It is evident.
\item [\ref{noco2}$\Rightarrow$\ref{noco3}:] It follows by Proposition \ref{omegprop}\ref{omegprop2}.
\item [\ref{noco3}$\Rightarrow$\ref{noco4}:] It is evident.
\item [\ref{noco4}$\Rightarrow$\ref{noco5}:] Let $x$ and $y$ be two pairwise elements in $A$. Assume by absurdum that $x^{\perp}\veebar y^{\perp}\nsubseteq A$. So $x^{\perp}\veebar y^{\perp}\subseteq \mathfrak{m}$, for some maximal filter $\mathfrak{m}$ of $\mathfrak{A}$. Applying Proposition \ref{1mineq}\ref{1mineq3}, it verifies that $x,y\notin D(\mathfrak{m})$; a contradiction.
\item [\ref{noco5}$\Rightarrow$\ref{noco6}:] It follows by Proposition \ref{compropo}.
\item [\ref{noco6}$\Rightarrow$\ref{noco7}:] Let $a\in (x\vee y)^{\perp}$. Let $b=a\vee x$. Obviously, $b$ and $y$ are pairwise. There exists $s\in A$ such that $s\in b^{\perp}$ and $\neg s\in y^{\perp}$. By \ref{res2} follows that $a\geq (a\vee s)\odot \neg s$. This establishes that $a\in x^{\perp}\veebar y^{\perp}$. The converse inclusion is evident.
\item [\ref{noco7}$\Rightarrow$\ref{noco8}:] It is evident.
\item [\ref{noco8}$\Rightarrow$\ref{noco1}:] Let $\mathfrak{m}$ and $\mathfrak{n}$ be distinct minimal prime filters of $\mathfrak{A}$. Consider $x\in \mathfrak{m}\setminus \mathfrak{n}$ and $y\in \mathfrak{n}\setminus \mathfrak{m}$. Using Proposition \ref{1mineq}\ref{1mineq3}, there exists $z\in x^{\perp}\setminus \mathfrak{m}$. Let $a=y\vee z$. So $(a\vee x)^{\perp}=A$, and this implies that $A=a^{\perp}\veebar x^{\perp}\subseteq \mathfrak{m}\veebar \mathfrak{n}$.
\end{proof}

\begin{theorem}\label{osublatfil}
Let $\mathfrak{A}$ be a residuated lattice. The following assertions are equivalent:
\begin{enumerate}
\item  [$(1)$ \namedlabel{osublatfil1}{$(1)$}] for any $F,G\in \Omega(\mathfrak{A})$, $F\vee^{\omega} G=A$ implies $F\veebar G=A$;
\item  [$(2)$ \namedlabel{osublatfil2}{$(2)$}] $\mathfrak{A}$ is mp;
\item  [$(3)$ \namedlabel{osublatfil3}{$(3)$}] for any $\mathcal{F}\subseteq \Omega(\mathfrak{A})$, $\veebar \mathcal{F}\in \Omega(\mathfrak{A})$;
%$\Omega(\mathfrak{A})$ is a complete sublattice of $\mathscr{F}(\mathfrak{A})$;
%for any family $\{I_i\}_{i\in I}$ of lattice ideals, $\veebar_{i\in I}\omega(I_i)=\omega(\curlyveeI_i)$;
\item  [$(4)$ \namedlabel{osublatfil4}{$(4)$}] $(\Omega(\mathfrak{A});\cap,\veebar)$ is a frame;
\item  [$(5)$ \namedlabel{osublatfil5}{$(5)$}] $(\gamma(\mathfrak{A});\cap,\veebar)$ is a lattice.
\end{enumerate}
\end{theorem}
\begin{proof}
\item  [\ref{osublatfil1}$\Rightarrow$\ref{osublatfil2}:] Let $x\vee y=1$ for some $x,y\in A$. Since $\gamma(\mathfrak{A})$ is a sublattice of $\Omega(\mathfrak{A})$ so we have $x^{\perp}\vee^{\omega} y^{\perp}=x^{\perp}\vee^{\Gamma} y^{\perp}=(x\vee y)^{\perp}=A$.
\item  [\ref{osublatfil2}$\Rightarrow$\ref{osublatfil3}:] Let $\{F_i\}_{i\in I}$ be a family of $\omega$-filters of $\mathfrak{A}$. Obviously, we have $\veebar_{i\in I} F_i\subseteq \omega(\curlyvee_{i\in I}I_i)$. Consider $a\in \omega(\curlyvee_{i\in I}I_i)$. Hence, there exists $x\in \curlyvee_{i\in I}I_i$ such that $a\in x^{\perp}$. This states that $x\leq x_{i_1}\vee\cdots\vee x_{i_n}$, for some integer $n$ and $x_{i_j}\in I_{i_j}$. We have the following sequence of formulas:
     \[
     \begin{array}{ll}
        x^{\perp}&\subseteq (x_{i_1}\vee\cdots\vee x_{i_n})^{\perp} \\
        &=x^{\perp}_{i_1}\veebar\cdots\veebar x^{\perp}_{i_n} \\
        &\subseteq  F_{i_1}\veebar\cdots\veebar F_{i_n}\\
        &\subseteq \veebar_{i\in I} F_i.
     \end{array}
     \]
\item  []\ref{osublatfil3}$\Rightarrow$\ref{osublatfil4}: It is evident.
\item  []\ref{osublatfil4}$\Rightarrow$\ref{osublatfil5}: It follows by Proposition \ref{omegprop}\ref{omegprop1}.
\item  []\ref{osublatfil5}$\Rightarrow$\ref{osublatfil1}: Let $F,G\in \Omega(\mathfrak{A})$ such that $F\vee^{\omega} G=A$. Since $\omega(I_F\curlyvee I_G)=A$, so $1\in I_F\curlyvee I_G$. This establishes that $f\vee g=1$, for some $f\in I_F$ and $g\in I_G$. Hence, $A=(f\vee g)^{\perp}=f^{\perp}\vee^{\Gamma} g^{\perp}=f^{\perp}\veebar g^{\perp}\subseteq F\veebar G$.
\end{proof}

%For a prime filter $P$ of an mp-residuated lattice $\mathfrak{A}$ the unique minimal prime filter of $\mathfrak{A}$ contained in $P$ denoted by $\mathfrak{m}_{P}$.
%\begin{theorem}\citep{rasouli2018n}\label{noco}
%Let $\mathfrak{A}$ be a residuated lattice. The following assertions are equivalent:
%\begin{enumerate}
%\item  [$(1)$ \namedlabel{noco1}{$(1)$}] $\mathfrak{A}$ is mp;
%\item  [$(2)$ \namedlabel{noco2}{$(2)$}] Any two distinct minimal prime filters of $\mathfrak{A}$ are comaximal;
%\item  [$(3)$ \namedlabel{noco3}{$(3)$}] for any $x,y\in A$, $x\vee y=1$ implies $x^{\perp}$ and $y^{\perp}$ are comaximal;
%\item  [$(4)$ \namedlabel{noco4}{$(4)$}] $(\Omega(\mathfrak{A});\cap,\veebar)$ is a frame;
%\item  [$(5)$ \namedlabel{noco5}{$(5)$}] $(\gamma(\mathfrak{A});\cap,\veebar)$ is a lattice.;
%\item  [$(6)$ \namedlabel{noco6}{$(6)$}] $\mathfrak{D(A)}\subseteq Min(\mathfrak{A})$.
%\end{enumerate}
%\end{theorem}

\citet[Proposition 2.1]{matlis1983minimal} gave a criterion for a ring to be PF and showed that a unitary commutative ring $\mathfrak{A}$ is PF if and only if for any maximal ideal $\mathfrak{m}$ of $\mathfrak{A}$, $\mathfrak{A_m}$ be an integral domain. Motivated by this, the following theorem, which is an immediate consequence of Proposition \ref{domchar} and Theorem \ref{noco}, can be extracted for mp-residuated lattices.
\begin{theorem}\label{matlchar}(Matlis's characterization)
Let $\mathfrak{A}$ be a residuated lattice. The following assertions are equivalent:
\begin{enumerate}
   \item [(1) \namedlabel{matlchar1}{(1)}] $\mathfrak{A}$ is mp;
   \item [(2) \namedlabel{matlchar2}{(2)}] $\mathfrak{A}/D(\mathfrak{p})$ is a domain, for any prime filter $\mathfrak{p}$ of $\mathfrak{A}$;
   \item [(3) \namedlabel{matlchar3}{(3)}] $\mathfrak{A}/D(\mathfrak{m})$ is a domain, for any maximal filter $\mathfrak{m}$ of $\mathfrak{A}$.
 \end{enumerate}
\end{theorem}

The next theorem gives some necessary and sufficient conditions for the collection of minimal prime filters in a residuated lattice to be a Hausdorff space with the dual hull-kernel topology.
\begin{theorem}\label{mpmpropd}
Let $\mathfrak{A}$ be a residuated lattice.  The following assertions are equivalent:
 \begin{enumerate}
   \item [(1) \namedlabel{mpmpropd1}{(1)}] $\mathfrak{A}$ is mp;
   \item [(2) \namedlabel{mpmpropd2}{(2)}] $Min_{d}(\mathfrak{A})$ is Hausdorff.
 \end{enumerate}
\end{theorem}
\begin{proof}
  \item [\ref{mpmpropd1}$\Rightarrow$\ref{mpmpropd2}:] Let $\mathfrak{m}$ and $\mathfrak{n}$ be two distinct minimal prime filters of $\mathfrak{A}$. So there exist $x\in \mathfrak{m}$ and $y\in \mathfrak{n}$ such that $x\odot y=0$. This follows that $h(x)\cap h(y)=\emptyset$, and this holds the result.
  \item [\ref{mpmpropd2}$\Rightarrow$\ref{mpmpropd1}:] Let $\mathfrak{m}$ and $\mathfrak{n}$ be two distinct minimal prime filters of $\mathfrak{A}$. So there exist $x,y\in A$ such that $\mathfrak{m}\in h(x)$, $\mathfrak{n}\in h(y)$, and $h(x\odot y)=\emptyset$. This shows that $A=x^{\perp\perp}\veebar y^{\perp\perp}\subseteq \mathfrak{m}\veebar \mathfrak{n}$. This holds the result.
\end{proof}
\begin{remark}\label{mpminhau}
By Proposition \ref{dhulkerinstr} and Theorem \ref{mpmpropd}, $\mathfrak{A}$ is an mp-residuated lattice if and only if $Min_{d}(\mathfrak{A})$ is a $T_{4}$ space.
\end{remark}
\begin{theorem}\label{minhausdclossp}
Let $\mathfrak{A}$ be a residuated lattice.  The following assertions are equivalent:
 \begin{enumerate}
   \item [(1) \namedlabel{minhausdclossp1}{(1)}] $\mathfrak{A}$ is mp;
   \item [(2) \namedlabel{minhausdclossp2}{(2)}] $h(\mathfrak{m})$ is closed in $Spec_{d}(\mathfrak{A})$, for any $\mathfrak{m}\in Min(\mathfrak{A})$.
 \end{enumerate}
\end{theorem}
\begin{proof}
  \item [\ref{minhausdclossp1}$\Rightarrow$\ref{minhausdclossp2}:] It follows by Proposition \ref{1mineq}\ref{1mineq2} and Theorem \ref{closefalzai}.

  %Let $\mathfrak{m}$ be a minimal prime filter of $\mathfrak{A}$. Consider $P\in d(\mathfrak{m})$. Assume that $\mathfrak{n}$ is a minimal prime filter contained in $P$. So there exist some elements $x\in \mathfrak{m}$ and $y\in \mathfrak{n}$ such that $x\odot y=0$. Thus $h(x)\cap h(y)=\emptyset$, and this implies that $P\in h(y)\subseteq d(x)\subseteq d(\mathfrak{m})$. This holds the result.
  \item [\ref{minhausdclossp2}$\Rightarrow$\ref{minhausdclossp1}:] Assume by absurdum that there exist two distinct minimal prime filters $\mathfrak{m}$ and $\mathfrak{n}$ of $\mathfrak{A}$ such that $\mathfrak{m}\veebar \mathfrak{n}\neq A$. This implies that there exists a prime filter $P$ containing in $\mathfrak{m}$ and $\mathfrak{n}$, and so $h(\mathfrak{m})\cap h(\mathfrak{n})\neq \emptyset$.
\end{proof}

Recall that a \textit{retraction} is a continuous mapping from a topological space into a subspace which preserves the position of all points in that subspace.
\begin{theorem}\label{mpculpro}
  Let $\mathfrak{A}$ be a residuated lattice. The following assertions are equivalent:
\begin{enumerate}
\item  [$(1)$ \namedlabel{mpculpro1}{$(1)$}] $\mathfrak{A}$ is mp;
\item  [$(2)$ \namedlabel{mpculpro2}{$(2)$}] $Min_{d}(\mathfrak{A})$ is a retraction of $Spec_{d}(\mathfrak{A})$.
\end{enumerate}
\end{theorem}
\begin{proof}
\item [\ref{mpculpro1}$\Rightarrow$\ref{mpculpro2}:] Define $f:Spec(\mathfrak{A})\longrightarrow Min(\mathfrak{A})$ by $f(\mathfrak{p})=\mathfrak{m}_{\mathfrak{p}}$. Set $H=\{\mathfrak{p}\in Spec(\mathfrak{A})\mid a\notin f(\mathfrak{p})\}$ and $X=(\bigcup H)^{c}$. Consider $a\in A$. Let $\mathfrak{p}\in f^{-1}(d_{m}(a))$. This implies that $\mathfrak{p}\in H$, and so $\mathfrak{p}\cap X=\emptyset$. Conversely, suppose that $\mathfrak{p}\cap X=\emptyset$. Let $a\in f(\mathfrak{p})$. So for any $\mathfrak{n}\in d_{m}(a)$ there exist $x_{\mathfrak{n}}\in \mathfrak{n}$ and $y_{\mathfrak{n}}\in f(\mathfrak{p})$ such that $x_{\mathfrak{n}}\odot y_{\mathfrak{n}}=0$. Obviously, $d_{m}(a)\subseteq \bigcup_{\mathfrak{n}\in d_{m}(a)} h(x_{\mathfrak{n}})$. Since $d_{m}(a)$ is a compact subspace of $Min_{d}(\mathfrak{A})$, so $d_{m}(a)\subseteq \bigcup_{\mathfrak{n}\in \Im} h(x_{\mathfrak{n}})=h(\bigvee_{\mathfrak{n}\in \Im} x_{\mathfrak{n}})$, where $\Im$ is a finite subset of $d_{m}(a)$. Letting $y=\bigodot_{\mathfrak{n}\in \Im} y_{\mathfrak{n}}$, we have $y\in f(\mathfrak{p})$ and $x\odot y=0$. So there exists a prime filter $Q$ of $\mathfrak{A}$ such that $Q\cap X=\emptyset$ and $y\in Q$. Since $a\notin f(Q)$, so $x\in Q$; a contradiction. This shows that $f^{-1}(d_{m}(a))=\{\mathfrak{p}\mid \mathfrak{p}\cap X=\emptyset\}$. So the result holds by \textsc{Remark} \ref{closd}.
\item [\ref{mpculpro2}$\Rightarrow$\ref{mpculpro1}:] Let $f:Spec_{d}(\mathfrak{A})\longrightarrow Min_{d}(\mathfrak{A})$ be a retraction and $\mathfrak{m}\in Min(\mathfrak{A})$.  Suppose that $\mathfrak{m}\subseteq \mathfrak{p}$, for some $\mathfrak{p}\in Spec(\mathfrak{A})$. By Lemma \ref{retractlemma}, we have $\mathfrak{m}\in cl_{d}^{Spec(\mathfrak{A})}(\mathfrak{p})$ and by continuity of $f$ and $T_1$ we obtain that
\[\mathfrak{m}=f(\mathfrak{m})\in f(cl_{d}^{Spec(\mathfrak{A})}(\mathfrak{p}))\subseteq cl_{d}^{Min(\mathfrak{A})}(f(\mathfrak{p}))=\{f(\mathfrak{p})\}.\]
This shows that $\mathfrak{m}$ is the unique minimal prime filter of $\mathfrak{A}$ contained in $\mathfrak{p}$.
\end{proof}
\begin{remark}\label{remardkretr}
By Theorem \ref{mpculpro}, if $\mathfrak{A}$ is an mp-residuated lattice, the map $Spec(\mathfrak{A})\rightsquigarrow Min(\mathfrak{A})$, which sends any prime filter $\mathfrak{p}$ of $\mathfrak{A}$ to the unique minimal prime filter of $\mathfrak{A}$ containing in it, is the unique retraction from $Spec_{d}(\mathfrak{A})$ into $Min_{d}(\mathfrak{A})$.
\end{remark}

The next result, which can be compared with Proposition \ref{minpro}, characterizes the clopen subsets of $Min_{d}(\mathfrak{A})$ where $\mathfrak{A}$ is an mp-residuated lattice.
\begin{corollary}\label{minmpproclo}
Let $\mathfrak{A}$ be an mp-residuated lattice. We have:
\[Clop(Min_{d}(\mathfrak{A}))=\{h(e)\cap Min(\mathfrak{A})\mid e\in \beta(\mathfrak{A})\}.\]
\end{corollary}
\begin{proof}
  By Theorem \ref{mpculpro}, there exists a retraction $f:Spec_{d}(\mathfrak{A})\longrightarrow Min_{d}(\mathfrak{A})$. Let $U\in Clop(Min_{d}(\mathfrak{A}))$. So $f^{\leftarrow}(U)\in Clop(Spec_{d}(\mathfrak{A}))$. Thus $f^{\leftarrow}(U)=h(e)$, for some $e\in \beta(\mathfrak{A})$, due to Proposition \ref{minpro}. This implies that $U=f^{\leftarrow}(U)\cap Min(\mathfrak{A})=h(e)\cap Min(\mathfrak{A})$. The converse is evident.
\end{proof}
\begin{theorem}\label{mp1culpro}
  Let $\mathfrak{A}$ be a residuated lattice. The following assertions are equivalent:
\begin{enumerate}
\item  [$(1)$ \namedlabel{mp1culpro1}{$(1)$}] $\mathfrak{A}$ is mp;
\item  [$(2)$ \namedlabel{mp1culpro2}{$(2)$}] $Spec_{d}(\mathfrak{A})$ is a normal space.
\end{enumerate}
\end{theorem}
\begin{proof}
\item [\ref{mp1culpro1}$\Rightarrow$\ref{mp1culpro2}:] Using Theorem \ref{mpculpro} and \textsc{Remark} \ref{remardkretr}, there exists a retraction $f:Spec_{d}(\mathfrak{A})\longrightarrow Min_{d}(\mathfrak{A})$, which sends any prime filter of $\mathfrak{A}$ to the unique minimal prime filter of $\mathfrak{A}$ contained in $\mathfrak{p}$, for any prime filter $\mathfrak{p}$ of $\mathfrak{A}$. By \textsc{Remark} \ref{mpminhau}, $Min(\mathfrak{A})$ is a $T_4$ space, and so $f$ is a closed map. Let $C_1$ and $C_2$ be two disjoint closed sets in $Spec_{d}(\mathfrak{A})$, so $f(C_1)$ and $f(C_2)$ are disjoint closed sets in $Min_{d}(\mathfrak{A})$. Since $Min_{d}(\mathfrak{A})$ is normal, there exist disjoint open neighbourhoods $N_1$ and $N_2$ of $f(C_1)$  and $f(C_2)$ in $Min_{d}(\mathfrak{A})$, respectively. One can see that  $f^{-1}(N_1)$ and $f^{-1}(N_2)$ are disjoint open neighbourhoods of $C_1$ and $C_2$ in $Spec_{d}(\mathfrak{A})$, respectively.
\item [\ref{mpculpro2}$\Rightarrow$\ref{mpculpro1}:] Let $\mathfrak{m}\in Min(\mathfrak{A})$. If $\mathfrak{p}\in Cl_{d}^{Spec(\mathfrak{A})}(\mathfrak{m})$, $\mathfrak{p}\subseteq \mathfrak{m}$, and this yields that $\mathfrak{p}=\mathfrak{m}$. This shows that $\{\mathfrak{m}\}$ is a closed subset of $Spec_{d}(\mathfrak{A})$. Now, let $\mathfrak{m}_{1},\mathfrak{m}_{2}\in Min(\mathfrak{A})$. Thus There exist $a,b\in A$ such that $h(a)$ and $h(b)$ are disjoint neighborhood of $\mathfrak{m}_{1}$ and $\mathfrak{m}_{2}$ in $Spec_{d}(\mathfrak{A})$, respectively. This shows that $h_{m}(a)$ and $h_{m}(b)$ are disjoint neighborhood of $\mathfrak{m}_{1}$ and $\mathfrak{m}_{2}$ in $Min_{d}(\mathfrak{A})$, respectively. This holds the result.
\end{proof}

Let $\mathfrak{A}$ be a residuated lattice. Consider the following relation $\imath=\{(\mathfrak{p},\mathfrak{q})\in X^{2}\mid \mathfrak{p}\veebar \mathfrak{q}\neq A\}$ on $X=Spec(\mathfrak{A})$. Obviously, $\imath$ is reflexive and symmetric. Let $\overline{\imath}$ be the transitive closure of $\imath$.
\begin{theorem}\label{31hausnorm}
  Let $\mathfrak{A}$ be a residuated lattice. The following assertions are equivalent:
  \begin{enumerate}
   \item [(1) \namedlabel{31hausnorm1}{(1)}] $\mathfrak{A}$ is mp;
   \item [(2) \namedlabel{31hausnorm2}{(2)}] for a given minimal prime filter $\mathfrak{m}$ of $\mathfrak{A}$, $\overline{\imath}(\mathfrak{m})=h(\mathfrak{m})$.
  % \item [(3) \namedlabel{31hausnorm3}{(3)}] if $\mathfrak{m}$ and $\mathfrak{n}$ are two distinct minimal prime filters of $\mathfrak{A}$, then $\Lambda(M\cap N)=\emptyset$.
 \end{enumerate}
\end{theorem}
\begin{proof}
  \item [\ref{31hausnorm1}$\Rightarrow$\ref{31hausnorm2}:] Let $\mathfrak{m}$ be a minimal prime filter of $\mathfrak{A}$. Consider $\mathfrak{p}\in \overline{\imath}(\mathfrak{m})$. So there exists a finite set $\{\mathfrak{p}_1,\cdots,\mathfrak{p}_n\}$ of elements of $Spec(\mathfrak{A})$ with $n\geq 2$ such that $\mathfrak{p}_1=\mathfrak{p}$, $\mathfrak{p}_n=\mathfrak{m}$, and $(\mathfrak{p}_{i},\mathfrak{p}_{i+1})\in \imath$, for all $1\leq i\leq n-1$. If $n=2$, then $\mathfrak{p}\veebar \mathfrak{m}\neq A$, and so $\mathfrak{m}\subseteq \mathfrak{p}$. Assume that $n>2$. We have $\mathfrak{p}_{n-2}\veebar \mathfrak{p}_{n-1}\neq A$ and $\mathfrak{m}\subseteq \mathfrak{p}_{n-1}$. This verifies that $(\mathfrak{p}_{n-2},\mathfrak{m})\in \imath$. Hence, in the equivalency $(\mathfrak{p},\mathfrak{m})\in \overline{\imath}$, the number of the involved primes is reduced to $n-1$. Therefore by the induction hypothesis, $\mathfrak{m}\subseteq \mathfrak{p}$. This shows that $\overline{\imath}(\mathfrak{m})\subseteq h(\mathfrak{m})$. The inverse inclusion is evident.
  \item [\ref{3hausnorm2}$\Rightarrow$\ref{3hausnorm1}:] It is evident.
  %\item [\ref{3hausnorm3}$\Rightarrow$\ref{3hausnorm1}:] Let $P$ be a prime filter of $\mathfrak{A}$ containing $Rad(\mathfrak{A})$ such that $P\subseteq M\cap N$, for some maximal filters $M$ and $N$ of $\mathfrak{A}$. It follows that $P\in \Lambda(M)\cap \Lambda(N)$ and so $M=N$. This holds the result.
\end{proof}

Let $A_{\tau}$ be a topological space, and $E$ be an equivalence relation on $A$. In the following, by $A_{\tau}/E$ we mean the quotient of the space $A_{\tau}$ modulo $E$. By \citet[p.90]{engelking1989general}, the quotient map $\pi:A_{\tau}\longrightarrow A_{\tau}/E$ is continuous, and a mapping $f$ of the quotient space $A_{\tau}/E$ to a topological space $B_{\zeta}$ is continuous if and only if the composition $f\circ \pi$ is continuous.
\begin{corollary}\label{minihomeo}
  Let $\mathfrak{A}$ be a residuated lattice. $\mathfrak{A}$ is mp if and only if the map $\eta:Min_{d}(\mathfrak{A})\longrightarrow Spec_{d}(\mathfrak{A})/\overline{\imath}$, given by $\mathfrak{m}\rightsquigarrow \overline{\imath}(\mathfrak{m})$, is a homeomorphism.
\end{corollary}
\begin{proof}
  Let $Min_{d}(\mathfrak{A})$ is a Hausdorff space. It is evident that $Spec_{d}(\mathfrak{A})/\overline{\imath}=\{\overline{\imath}(\mathfrak{m})\mid \mathfrak{m}\in Min(\mathfrak{A})\}$, and this implies that $\eta$ is a surjection. The injectivity of $\eta$ follows by Theorem \ref{31hausnorm}, and the continuity of it follows by $\eta=\pi\circ i$, where $i$ is the inclusion map. By \textsc{Remark} \ref{remardkretr} and Theorem \ref{mpculpro} follows that $\eta^{-1}\circ \pi$ is a retraction, and this verifies the continuity of $\eta^{-1}$, see \citet[Proposition 4.2.4]{engelking1989general}. This shows that $\eta$ is a homeomorphism. Conversely, let $\eta:Min_{d}(\mathfrak{A})\longrightarrow Spec_{d}(\mathfrak{A})/\overline{\imath}$ be a homeomorphism. Obviously, $\eta^{-1}\circ \pi$ is a retraction, and so $Min_{d}(\mathfrak{A})$ is a Hausdorff space due to Theorem \ref{mpculpro}.
\end{proof}

Let $\mathfrak{A}$ be a residuated lattice. Consider the relation $\jmath=\{(\mathfrak{p},\mathfrak{q})\in X^{2}\mid \dot{\mathfrak{p}}\curlyvee \dot{\mathfrak{q}}\neq A\}$ on $X=Spec(\mathfrak{A})$. Obviously, $\jmath$ is reflexive and symmetric. Let $\overline{\jmath}$ be the transitive closure of $\jmath$.
\begin{remark}
  For prime filters $\mathfrak{p}$ and $\mathfrak{q}$ of a residuated lattice $\mathfrak{A}$. One can see that, using \cite[Proposition 3.5]{rasouli2018n}, $\dot{\mathfrak{p}}\curlyvee \dot{\mathfrak{q}}=A$ if and only if $D(\mathfrak{p})\veebar D(\mathfrak{q})=A$.
\end{remark}
\begin{theorem}\label{3hausnorm}
  Let $\mathfrak{A}$ be a residuated lattice. The following assertions are equivalent:
  \begin{enumerate}
   \item [(1) \namedlabel{3hausnorm1}{(1)}] $\mathfrak{A}$ is mp;
   \item [(2) \namedlabel{3hausnorm2}{(2)}] for a given minimal prime filter $\mathfrak{m}$ of $\mathfrak{A}$, $\overline{\jmath}(\mathfrak{m})=h(\mathfrak{m})$.
   %\item [(3) \namedlabel{3hausnorm3}{(3)}] if $\mathfrak{m}$ and $\mathfrak{n}$ are two distinct minimal prime filters of $\mathfrak{A}$, then $h(\mathfrak{m}) \cap h(\mathfrak{n})=\emptyset$.
 \end{enumerate}
\end{theorem}
\begin{proof}
  \item [\ref{3hausnorm1}$\Rightarrow$\ref{3hausnorm2}:] Let $\mathfrak{m}$ be a minimal prime filter of $\mathfrak{A}$. Consider $\mathfrak{p}\in \overline{\jmath}(\mathfrak{m})$. So there exists a finite set $\{\mathfrak{p}_1,\cdots,\mathfrak{p}_n\}$ of elements of $Spec(\mathfrak{A})$ with $n\geq 2$ such that $\mathfrak{p}_1=P$, $\mathfrak{p}_n=\mathfrak{m}$, and $(\mathfrak{p}_{i},\mathfrak{p}_{i+1})\in \jmath$, for all $1\leq i\leq n-1$. If $n=2$, then $\dot{\mathfrak{p}}\curlyvee \mathfrak{\dot{m}}\neq A$, and so $\mathfrak{m}\subseteq \mathfrak{p}$ due to Proposition \ref{1mineq}. Assume that $n>2$. We have $\dot{\mathfrak{p}_{n-2}}\curlyvee \dot{\mathfrak{p}_{n-1}}\neq A$ and $\mathfrak{m}\subseteq \mathfrak{p}_{n-1}$. Using Zorn's lemma, it verifies that $\dot{\mathfrak{p}_{n-2}}\curlyvee \dot{\mathfrak{p}_{n-1}}\subseteq \mathfrak{c}$, for a maximal $\vee$-closed set of $\mathfrak{A}$. Applying Proposition \ref{1mineq} and the hypothesis, it shows that $\mathfrak{m}=\dot{\mathfrak{c}}$. This verifies that $(\mathfrak{p}_{n-2},\mathfrak{m})\in \jmath$. Hence, in the equivalency $(\mathfrak{p},\mathfrak{m})\in \overline{\jmath}$, the number of the involved primes is reduced to
$n-1$. Therefore by the induction hypothesis, $\mathfrak{m}\subseteq \mathfrak{p}$. This shows that $\overline{\jmath}(\mathfrak{m})\subseteq h(\mathfrak{m})$. The inverse inclusion is evident.
  \item [\ref{3hausnorm2}$\Rightarrow$\ref{3hausnorm1}:] It is evident.
 % \item [\ref{3hausnorm3}$\Rightarrow$\ref{3hausnorm1}:] Let $P$ be a prime filter of $\mathfrak{A}$ containing $\mathfrak{m}$ and $\mathfrak{n}$. It follows that $P\in \mathscr{S}(\mathfrak{m}\veebar \mathfrak{n})$, and so $\mathfrak{m}=\mathfrak{n}$. This holds the result.
\end{proof}

The proof of the following corollary is analogous to the proof of Corollary \ref{minihomeo}, and so it is left to the reader.
\begin{corollary}\label{minjhomeo}
  Let $\mathfrak{A}$ be a residuated lattice. $\mathfrak{A}$ is mp if and only if the map $\eta:Min_{d}(\mathfrak{A})\longrightarrow Spec_{d}(\mathfrak{A})/\overline{\jmath}$, given by $\mathfrak{m}\rightsquigarrow \overline{\jmath}(\mathfrak{m})$, is a homeomorphism.
\end{corollary}

%%%%%%%%%%%%%%%%%%%%%%%%%%%%%%%%%%%%%%%%%%%%%%%%%%%%%%%%%%%%%%%%%%%%%%%%%%%%%%%%%%%%%%%%%
\section{The pure spectrum of an mp-residuated lattice}\label{sec4}

This section deals with the pure spectrum of an mp-residuated lattice.  For the basic facts concerning pure filters of a residuated lattice, the reader is referred to \cite{rasouli2021rickart}.

For any filter $F$ of a residuated lattice $\mathfrak{A}$, set $\sigma(F)=k\mathscr{G}h(F)$.
\begin{proposition}\label{unitpro}\cite[Propositions 5.2 \& 5.4]{rasouli2021rickart}
  Let $\mathfrak{A}$ be a residuated lattice. The following assertions hold:
\begin{enumerate}
\item  [$(1)$ \namedlabel{unitpro1}{$(1)$}] $\sigma(F)=\{a\in A\mid F\veebar a^{\perp}=A\}$, for any filter $F$ of $\mathfrak{A}$;
%\item  [$(2)$ \namedlabel{unitpro2}{$(2)$}] $\sigma(F)$ is a filter of $\mathfrak{A}$ contained in $F$, for any filter $F$ of $\mathfrak{A}$;
\item  [$(2)$ \namedlabel{unitpro3}{$(2)$}] $F\subseteq G$ implies $\sigma(F)\subseteq \sigma(G)$, for any filters $F$ and $G$ of $\mathfrak{A}$;
%\item  [$(4)$ \namedlabel{unitpro4}{$(4)$}] $\sigma(\mathfrak{p})\subseteq D(\mathfrak{p})$, for any prime filter $\mathfrak{p}$ of $\mathfrak{A}$;
\item  [$(3)$ \namedlabel{unitpro5}{$(3)$}] $\sigma(\mathfrak{m})=D(\mathfrak{m})$, for any maximal filter $\mathfrak{m}$ of $\mathfrak{A}$;
%\item  [$(6)$ \namedlabel{unitpro6}{$(6)$}] $\sigma(F\cap G)=\sigma(F)\cap \sigma(G)$, for any filters $F$ and $G$ of $\mathfrak{A}$;
%\item  [$(7)$ \namedlabel{unitpro7}{$(7)$}] $\veebar_{F\in \mathcal{F}}\sigma(F)\subseteq \sigma(\veebar \mathcal{F})$, for any family of filters of $\mathfrak{A}$.
\end{enumerate}
\end{proposition}

Let $\mathfrak{A}$ be a residuated lattice. A filter $F$ of $\mathfrak{A}$ is called \textit{pure} provided that $\sigma(F)=F$. The set of pure filters of $\mathfrak{A}$ is denoted by $\sigma(\mathfrak{A})$. It is obvious that $\{1\},A\in \sigma(\mathfrak{A})$.
\begin{proposition}\label{sigmnorpr}
  Let $\mathfrak{A}$ be an mp-residuated lattice and $F$ a filter of $\mathfrak{A}$. $\sigma(F)$ is a pure filter of $\mathfrak{A}$.
\end{proposition}
\begin{proof}
Let $x\in \sigma(F)$. Applying Proposition \ref{unitpro}\ref{unitpro1}, it follows that $F\veebar x^{\perp}=A$. So $f\odot y=0$, for some $f\in F$ and $y\in x^{\perp}$. By Proposition \ref{noco} there exist $a\in A$ such that $a\in x^{\perp}$ and $\neg a\in y^{\perp}$. This implies that $\neg a\in \sigma(\mathfrak{A})$, and so $x\in \sigma(\sigma(F))$.
\end{proof}

The following theorem gives some criteria for mp-residuated lattices by pure filters, inspired by the
one obtained for bounded distributive lattices by \citet[Theorem 2.11]{cornish1977ideals}.
%\begin{theorem}\label{norgammsig}
%  Let $\mathfrak{A}$ be a residuated lattice. The following assertions are equivalent:
%  \begin{enumerate}
%\item  [$(1)$ \namedlabel{norgammsig1}{$(1)$}] $\mathfrak{A}$ is mp;
%\item  [$(2)$ \namedlabel{norgammsig2}{$(2)$}] $\Omega(\mathfrak{A})\subseteq \sigma(\mathfrak{A})$;
%\item  [$(3)$ \namedlabel{norgammsig3}{$(3)$}] $\gamma(\mathfrak{A})\subseteq \sigma(\mathfrak{A})$.
%\end{enumerate}
%\end{theorem}
%\begin{proof}
%\item [\ref{norgammsig1}$\Rightarrow$\ref{norgammsig2}:] Let $F$ be an $\omega$-filter of $\mathfrak{A}$. So $F=\omega(I)$, for some ideal $I$ of $\ell(\mathfrak{A})$. Consider $x\in F$. So $x\in a^{\perp}$, for some $a\in I$. By Propositions \ref{omegprop}\ref{omegprop1} and \ref{noco}\ref{noco4} follows that $A=x^{\perp}\veebar a^{\perp}\subseteq x^{\perp}\veebar F$. This holds the result.
%\item [\ref{norgammsig2}$\Rightarrow$\ref{norgammsig3}:] By Propositions \ref{omegprop}\ref{omegprop2}, it is evident.
%\item [\ref{norgammsig3}$\Rightarrow$\ref{norgammsig1}:] Let $x\vee y=1$. So $x\in y^{\perp}=\sigma(y^{\perp})$ and this implies that $x^{\perp}\veebar y^{\perp}=A$. Hence the result holds by Proposition \ref{noco}\ref{noco4}..
%\end{proof}
 \begin{theorem}\label{norgammsig}
  Let $\mathfrak{A}$ be a residuated lattice. The following assertions are equivalent:
  \begin{enumerate}
\item  [$(1)$ \namedlabel{norgammsig1}{$(1)$}] $\mathfrak{A}$ is mp;
\item  [$(2)$ \namedlabel{norgammsig2}{$(2)$}] $\Omega(\mathfrak{A})\subseteq \sigma(\mathfrak{A})$;
\item  [$(3)$ \namedlabel{norgammsig3}{$(3)$}] $\gamma(\mathfrak{A})\subseteq \sigma(\mathfrak{A})$.
\end{enumerate}
\end{theorem}
\begin{proof}
\item [\ref{norgammsig1}$\Rightarrow$\ref{norgammsig2}:] Let $F$ be an $\omega$-filter of $\mathfrak{A}$. So $F=\omega(I)$, for some ideal $I$ of $\ell(\mathfrak{A})$. Consider $x\in F$. So $x\in a^{\perp}$, for some $a\in I$. By Propositions \ref{omegprop}\ref{omegprop1} and \ref{noco}\ref{noco4} follows that $A=x^{\perp}\veebar a^{\perp}\subseteq x^{\perp}\veebar F$. This holds the result.
\item [\ref{norgammsig2}$\Rightarrow$\ref{norgammsig3}:] By Propositions \ref{omegprop}\ref{omegprop2}, it is evident.
\item [\ref{norgammsig3}$\Rightarrow$\ref{norgammsig1}:] Let $x\vee y=1$. So $x\in y^{\perp}=\sigma(y^{\perp})$ and this implies that $x^{\perp}\veebar y^{\perp}=A$. Hence the result holds by Proposition \ref{noco}.
\end{proof}
\begin{remark}
  \citet[Theorem 1]{al1987some} showed that a unitary commutative ring is a PF ring if and only if any its annulet is a pure ideal. Thus, if we define PF-residuated lattices as those ones in which any coannulet is a pure filter, Theorem \ref{norgammsig} verifies that the class of PF residuated lattices coincides with the class of mp-residuated lattices.
\end{remark}
\begin{lemma}\label{comxpureprime}
  Let $\mathfrak{A}$ be a residuated lattice. Any two distinct elements of the set $Spec(\mathfrak{A})\cap \sigma(\mathfrak{A})$ are comaximal.
\end{lemma}
\begin{proof}
  Let $\mathfrak{p}_1$ and $\mathfrak{p}_2$ be two distinct elements of the set $Spec(\mathfrak{A})\cap \sigma(\mathfrak{A})$. Consider $x\in \mathfrak{p}_1\setminus \mathfrak{p}_2$. So $\mathfrak{p}_1\veebar x^{\perp}=A$ and $x^{\perp}\subseteq \mathfrak{p}_2$. This holds the result.
\end{proof}
\begin{theorem}\label{norgammsige}
  Let $\mathfrak{A}$ be a residuated lattice. The following assertions are equivalent:
  \begin{enumerate}
\item  [$(1)$ \namedlabel{norgammsige1}{$(1)$}] $\mathfrak{A}$ is mp;
\item  [$(2)$ \namedlabel{norgammsige2}{$(2)$}] $D(\mathfrak{p})$ is a pure filter of $\mathfrak{A}$, for any prime filter $\mathfrak{p}$ of $\mathfrak{A}$;
\item  [$(3)$ \namedlabel{norgammsige3}{$(3)$}] $D(\mathfrak{m})$ is a pure filter of $\mathfrak{A}$, for any maximal filter $\mathfrak{m}$ of $\mathfrak{A}$;
\item  [$(4)$ \namedlabel{norgammsige4}{$(4)$}] $Min(\mathfrak{A})\subseteq \sigma(\mathfrak{A})$.
\end{enumerate}
\end{theorem}
\begin{proof}
\item [\ref{norgammsige1}$\Rightarrow$\ref{norgammsige2}:] It follows by Theorem \ref{norgammsig}.
\item [\ref{norgammsige2}$\Rightarrow$\ref{norgammsige3}:] It is evident.
\item [\ref{norgammsige3}$\Rightarrow$\ref{norgammsige4}:] It follows, with a little bit of effort, by Proposition \ref{omegprop}\ref{omegprop3}.
\item [\ref{norgammsige4}$\Rightarrow$\ref{norgammsige1}:] It follows by Proposition \ref{noco} and Lemma \ref{comxpureprime}.
\end{proof}

%%%%%%%%%%%%%%%%%%%%%%%%%%%%%%%%%%%%%%%%%%%%%%%%%%%%%%%
 Let $\mathfrak{A}$ be a residuated lattice. Recall \citep{rasouli2021rickart} that a proper pure filter of $\mathfrak{A}$ is called \textit{purely-maximal} provided that it is a maximal element in the set of proper and pure filters of $\mathfrak{A}$. The set of purely-maximal filters of $\mathfrak{A}$ shall be denoted by $Max(\sigma(\mathfrak{A}))$. A proper pure filter $\mathfrak{p}$ of $\mathfrak{A}$ is called \textit{purely-prime} provided that $F_1\cap F_2\subseteq \mathfrak{p}$ implies $F_1\subseteq \mathfrak{p}$ or $F_2\subseteq \mathfrak{p}$, for any $F_1,F_2\in \sigma(\mathfrak{A})$. The set of all purely-prime filters of $\mathfrak{A}$ shall be denoted by $Spp(\mathfrak{A})$. It is obvious that $Max(\sigma(\mathfrak{A}))\subseteq Spp(\mathfrak{A})$. Zorn's lemma ensures that any proper pure filter is contained in a purely-maximal filter, and so in a purely-prime filter.
\begin{theorem}\label{normpurprimxa}
  Let $\mathfrak{A}$ be a residuated lattice. The following assertions are equivalent:
\begin{enumerate}
\item  [(1) \namedlabel{normpurprimxa1}{(1)}] $\mathfrak{A}$ is mp;
\item  [(2) \namedlabel{normpurprimxa2}{(2)}] $Min(\mathfrak{A})=Max(\sigma(\mathfrak{A}))$.
\end{enumerate}
\end{theorem}
\begin{proof}
\item [\ref{normpurprimxa1}$\Rightarrow$\ref{normpurprimxa2}:] Let $\mathfrak{m}$ be a minimal prime filter of $\mathfrak{A}$. By Theorem \ref{norgammsige} follows that $\mathfrak{m}$ is a pure filter of $\mathfrak{A}$. Thus there exists $\mathfrak{n}\in Max(\sigma(\mathfrak{A}))$ containing $\mathfrak{m}$. Let $a\in \mathfrak{n}$. So there exists $b\in a^{\perp}$ such that $\neg b\in \mathfrak{n}$. This implies that $b\notin \mathfrak{m}$, and so $a\in \mathfrak{m}$. Conversely, let $\mathfrak{p}$ be a purely-maximal filter of $\mathfrak{A}$. So $\mathfrak{p}\subseteq \mathfrak{n}$, for some $\mathfrak{n}\in Max(\mathfrak{A})$. Using Theorem \ref{noco}, Proposition \ref{unitpro}(\ref{unitpro3} \& \ref{unitpro5}), and Theorem \ref{norgammsige}, it shows that $\mathfrak{p}=D(\mathfrak{n})\in Min(\mathfrak{A})$.
\item [\ref{normpurprimxa2}$\Rightarrow$\ref{normpurprimxa1}:] It is evident by Theorem \ref{norgammsige}.
\end{proof}
%%%%%%%%%%%%%%%%%%%%%%%%%%%%%%%%%%%%%%%%%%%%%%%%%%%%%%%%%%%%%%%%%%%%%%%%%%%%%%%%%%%%%%%%
The following result generalized and improved \cite[Theorem 1.8]{al1989pure} to residuated lattices.
\begin{proposition}\label{pureinterd}
  Let $\mathfrak{A}$ be an mp-residuated lattice and $F$ a proper pure filter of $\mathfrak{A}$. We have
  \[F=k(Min(\mathfrak{A})\cap h(F)).\]
\end{proposition}
\begin{proof}
  By Theorem \ref{normpurprimxa}, $Min(\mathfrak{A})\cap h(F)\neq \emptyset$. Consider $a\in k(Min(\mathfrak{A})\cap h(F))$. Assume that $a^{\perp}\veebar F$ is proper. Thus $a^{\perp}\veebar F\subseteq \mathfrak{n}$, for some maximal filter $\mathfrak{n}$ of $\mathfrak{A}$. Let $\mathfrak{m}$ be a minimal prime filter of $\mathfrak{A}$ contained in $\mathfrak{n}$. This implies that $F\subseteq \mathfrak{m}$, and so $\neg b\in \mathfrak{n}$, for some $b\in a^{\perp}$; a contradiction.
\end{proof}
%%%%%%%%%%%%%%%%%%%%%%%%%%%%%%%%%%%%%%%%%%%%%%%%%%%%%%%%%%%%%%%%%%%%%%%%%%%%%%%%%%%%%%%%
The pure ideals of a PF ring are characterized in \citet[Theorems 2.4 and 2.5]{al1988pure}. These results have been improved and generalized to residuated lattices in Theorem \ref{mppurefcl} and Proposition \ref{mppure}.
\begin{theorem}\label{mppurefcl}
Let $\mathfrak{A}$ be an mp-residuated lattice. The pure filters of $\mathfrak{A}$ are precisely of the form $\bigcap_{\mathfrak{m}\in Min(\mathfrak{A})\cap \mathcal{C}}\mathfrak{m}$,where $\mathcal{C}$ runs over closed subsets of $Spec_{d}(\mathfrak{A})$.
\end{theorem}
\begin{proof}
let $a\in G:=\bigcap\{\mathfrak{m}\mid \mathfrak{m}\in Min(\mathfrak{A})\cap \mathcal{C}\}$, in which $\mathcal{C}$ is a closed subset of $Spec_{d}(\mathfrak{A})$. So for any $\mathfrak{m}\in Min(\mathfrak{A})\cap \mathcal{C}$, we have $\mathfrak{m}\veebar a^{\perp}=A$. By absurdum, assume that $G\veebar a^{\perp}\neq A$. So $G\veebar a^{\perp}$ is contained in a maximal filter $\mathfrak{n}$. Let $\mathfrak{o}$ be a minimal prime filter of $\mathfrak{A}$ contained in $\mathfrak{m}$. Obviously, $\mathfrak{o}\notin \mathcal{C}$. So for any $\mathfrak{m}\in  Min(\mathfrak{A})\cap \mathcal{C}$, there exist $x_{\mathfrak{m}}\in \mathfrak{m}$ and $y_{\mathfrak{m}}\in \mathfrak{o}$ such that $x_{\mathfrak{m}}\odot y_{\mathfrak{m}}=0$. Since $\mathcal{C}$ is stable under the generalization, so $\mathcal{C}\subseteq \bigcup_{\mathfrak{m}\in Min(\mathfrak{A})\cap \mathcal{C}} h(x_{\mathfrak{m}})$. By Proposition \ref{dhulkerinstr} follows that $\mathcal{C}$ is compact. So there exist a finite number $\mathfrak{m}_{1},\cdots,\mathfrak{m}_{n}\in Min(\mathfrak{A})\cap \mathcal{C}$ such that $\mathcal{C}\subseteq \bigcup_{i=1}^{n} h(x_{\mathfrak{m}_{i}})$. Set $x=\bigvee_{i=1}^{n} x_{\mathfrak{m}_{i}}$ and $y=\bigodot_{i=1}^{n}y_{\mathfrak{m}_{i}}$. Routinely, one can see that $0=x\odot y\in G\veebar \mathfrak{o}$; a contradiction. The converse follows by Proposition \ref{pureinterd}.
\end{proof}
%%%%%%%%%%%%%%%%%%%%%%%%%%%%%%%%%%%%%%%%%%%%%%%%%%%%%%%%%%%%%
Let $\mathfrak{A}$ be a residuated lattice. For any filter $F$ of $\mathfrak{A}$, we set
\[\rho(F)=\underline{\bigvee}\{G\in \sigma(\mathfrak{A})\mid G\subseteq F\},\]
and it is called \textit{the pure part of $F$}. Definitely, the pure part of a filter is the largest pure filter contained in it.

%The following proposition has a routine verification, and so its proof is left to the reader.
\begin{proposition}\label{rfilter}
  Let $\mathfrak{A}$ be a residuated lattice. Then
  \[\bigcap \{\rho(\mathfrak{m})\mid \mathfrak{m}\in Max(\mathfrak{A})\}=\{1\}.\]
\end{proposition}
\begin{proof}
It is an immediate consequence of \citet[Proposition 4.5(5)]{rasouli2022gelfand}.
\end{proof}
%\begin{proof}
%\ref{rfilter1} and \ref{rfilter2} follows by \citet[Proposition 4.14]{rasouli2021rickart}.
%\item [\ref{rfilter3}:] It is evident by Proposition \ref{sigmapro}\ref{sigmapro2}.
%\item [\ref{rfilter4}:] It is evident by Theorem \ref{sigmfiltlatt}\ref{sigmfiltlatt1}.
%\end{proof}
\begin{proposition}\label{mppureco1}
Let $\mathfrak{A}$ be an mp-residuated lattice and $a\in A$. Then
\begin{center}
  $a^{\perp}\cap F_{a}=\{1\}$, where $F_{a}=\bigcap_{\mathfrak{m}\in Max(\mathfrak{A})\cap h(a)}\rho(\mathfrak{m})$.
\end{center}
\end{proposition}
\begin{proof}
With a little bit of effort, it follows by Theorem \ref{norgammsig} and Proposition \ref{rfilter}.
\end{proof}
\begin{corollary}\label{mppu1re}
  If $\mathfrak{m}$ is a minimal prime filter of an mp-residuated lattice $\mathfrak{A}$, then $\mathfrak{m}=\underline{\bigvee}_{a\in \mathfrak{m}}F_{a}$.
\end{corollary}
\begin{proof}
Let $a\in \mathfrak{m}$. So $b\in a^{\perp}$, for some $b\notin \mathfrak{m}$. This implies that $a\in F_{\neg b}$. The reverse inclusion is deduced from Corollary \ref{mppureco1}.
\end{proof}
%%%%%%%%%%%%%%%%%%%%%%%%%%%%%%%%%%%%%%%%%%%%%%%%%%%%%%%%%%%%%%%%%%%%%%%%%%%%%%%%%%%%%%%%%%%%%%%%%%%%%%%%%%%%%%%%%%%%%%%%%%%%%
\begin{proposition}\label{mppure}
Let $\mathfrak{A}$ be an mp-residuated lattice. The pure filters of $\mathfrak{A}$ are precisely of the form $\bigcap_{\mathfrak{m}\in Max(\mathfrak{A})\cap h(F)}\rho(\mathfrak{m})$, where $F$ is a filter of $\mathfrak{A}$.
\end{proposition}
\begin{proof}
Let $\mathcal{C}=\{P\in Spec(\mathfrak{A})\mid P\cap \neg F=\emptyset\}$. One can see that $Max(\mathfrak{A})\cap h(F)=Min(\mathfrak{A})\cap \mathcal{C}$. This establishes the result due to \textsc{Remark} \ref{closd} and Theorem \ref{mppurefcl}.
\end{proof}

\citet[Theorem 3.5]{al1988pure} proved that every purely prime ideal of a PF ring is purely maximal. Now we provide an alternative proof to the following interesting result.
\begin{theorem}\label{mpminspp}
  Let $\mathfrak{A}$ be an mp residuated lattice. Then
  \[Spp(\mathfrak{A})\subseteq Max(\sigma(\mathfrak{A})).\]
\end{theorem}
\begin{proof}
Let $\mathfrak{p}$ be a purely prime filter of $\mathfrak{A}$. So $\mathfrak{p}\subseteq \mathfrak{m}$, for some $\mathfrak{m}\in Max(\sigma(\mathfrak{A}))$. By Theorem \ref{normpurprimxa} we have $\mathfrak{m}\in Min(\mathfrak{A})$. Let $a\in \mathfrak{m}$. By Proposition \ref{1mineq}\ref{1mineq3} we have $a^{\perp}\nsubseteq \mathfrak{m}$. By Proposition \ref{mppureco1} follows that $a^{\perp}\cap F_{a}\subseteq P$. By Theorem \ref{norgammsig} and Proposition \ref{mppure}, respectively, follows that $a^{\perp}$ and $F_{a}$ are pure filters. This implies that $F_{a}\subseteq \mathfrak{p}$. Hence by Corollary \ref{mppu1re} follows that  $\mathfrak{m}=\underline{\bigvee}_{a\in \mathfrak{m}}F_{a}\subseteq \mathfrak{p}$.
\end{proof}

The following theorem is a direct consequence of Theorems \ref {norgammsige}, \ref {normpurprimxa} and \ref {mpminspp}. So its proof is left to the reader.
\begin{theorem}\label{mp2minspp}
  Let $\mathfrak{A}$ be a residuated lattice. The following assertions are equivalent:
\begin{enumerate}
\item  [(1) \namedlabel{mp2minspp1}{(1)}] $\mathfrak{A}$ is mp;
\item  [(2) \namedlabel{mp2minspp2}{(2)}] $Min(\mathfrak{A})=Spp(\mathfrak{A})$.
\end{enumerate}
\end{theorem}
%\begin{proof}
%\item [\ref{mp2minspp1}$\Rightarrow$\ref{mp2minspp2}:] It follows by Theorems \ref{normpurprimxa} \& \ref{mpminspp}.
%\item [\ref{mp2minspp2}$\Rightarrow$\ref{mp2minspp1}:] It is evident by Theorem \ref{norgammsige}.
%\end{proof}
%%%%%%%%%%%%%%%%%%%%%%%%%%%%%%%%%%%%%%%%%%%%%%%%%%%%%%%%%%%%%%%%%%%%%%%%%%%%%%%%%%%%%%%%%%%%%%%%%%%%%%%%%%%%%%%%%%%%%%%%%%%%%

For each pure filter $F$ of $\mathfrak{A}$ we set $d_{p}(F)=\{P\in Spp(\mathfrak{A})\mid F\nsubseteq P\}$. $Spp(\mathfrak{A})$ can be topologized by taking the set $\{d_{p}(F)\mid F\in \sigma(\mathfrak{A})\}$ as the open sets. The set $Spp(\mathfrak{A})$ endowed with this topology is called the \textit{pure spectrum} of $\mathfrak{A}$. It is obvious that the closed subsets of the pure spectrum are precisely of the form $h_{p}(F) =\{P\in Spp(\mathfrak{A})\mid F\subseteq P\}$, where $F$ runs over pure filters of $\mathfrak{A}$.

The next result, which can be compared with Proposition \ref{dhulkerinstr}, shows that the pure spectrum of a residuated lattice is a compact space.
\begin{theorem}\cite[Theorem 4.22]{rasouli2021rickart}\label{purespeccomp}
  Let $\mathfrak{A}$ be a residuated lattice. $Spp(\mathfrak{A})$ is a compact space.
\end{theorem}

The next theorem gives a criterion for a residuated lattice to be mp, inspired by the one obtained for unitary commutative rings by \citet[Theorem 5.4]{tarizadeh2020purely}.
\begin{theorem}\label{equmpflatmin}
  Let $\mathfrak{A}$ be a residuated lattice. The following assertions are equivalent:
\begin{enumerate}
\item  [(1) \namedlabel{equmpflatmin1}{(1)}] $\mathfrak{A}$ is mp;
\item  [(2) \namedlabel{equmpflatmin2}{(2)}] the identity map $\iota:Spp(\mathfrak{A})\longrightarrow Min_{d}(\mathfrak{A})$ is a homeomorphism.
\end{enumerate}
\end{theorem}
\begin{proof}
\item [\ref{equmpflatmin1}$\Rightarrow$\ref{equmpflatmin2}:] Consider the identity map $\iota:Spp(\mathfrak{A})\longrightarrow Min(\mathfrak{A})$. Using Theorem \ref{normpurprimxa}, it follows that $\iota$ is a well-defined bijection. One can see that $Min(\mathfrak{A})\cap h(a)=d_{p}(a^{\perp})$, for any $a\in A$, which implies that $\iota$ is continuous. By Theorems \ref{mpmpropd} and \ref{purespeccomp}, it follows that $Min_{d}(\mathfrak{A})$ is Hausdorff, and $Spp(\mathfrak{A})$ is compact, respectively. This holds the result due to \citet[Theorem 3.1.13]{engelking1989general}.
\item [\ref{equmpflatmin2}$\Rightarrow$\ref{equmpflatmin1}:] It is evident by Theorem \ref{normpurprimxa}.
\end{proof}

Using Theorem \ref{purespeccomp},  the pure spectrum of a residuated lattice is compact (not necessarily Hausdorff). The following result verifies that the pure spectrum of an mp-residuated lattice is Hausdorff.
\begin{corollary}\label{gelspphau}
  Let $\mathfrak{A}$ be an mp-residuated lattice. $Spp(\mathfrak{A})$ is a Hausdorff space.
\end{corollary}
\begin{proof}
It is an immediate consequence of Theorems \ref{mpmpropd} \& \ref{equmpflatmin}.
\end{proof}

%%%%%%%%%%%%%%%%%%%%%%%%%%%%%%%%%%%%%%%%%%%%%%%%%%%%%%%%%%%%%%%%%%%%%%%%%%%%%%%%%%%%%%%%

%%%%%%%%%%%%%%%%%%%%%%%%%%%%%%%%%%%%%%%%%%%%%%%%%%%%%%%%%%%%%%%%%%%%%%%%%%%%%%%%%%%%%%%%%%%%%%%%%%%%%%%%%%%%%%%%%%%%%%

\end{document}